\documentclass[12pt]{article}
\usepackage[active]{srcltx}
\usepackage[title]{appendix}
\usepackage{a4wide,amsfonts,amsmath,amsthm,amssymb,caption,pstricks,latexsym,delarray,empheq,hyperref,float,graphicx,setspace,subcaption,tabularx,url,enumitem, comment, tikz, natbib,environ,multirow, empheq, bbm, float,  titlesec,changepage,mathrsfs,chngcntr,geometry}

\newcommand{\R}{\mathbb{R}}
\newcommand{\probset}{\mathscr{P}}
\newcommand{\B}{\mathcal{B}}
\newcommand{\1}{\mathbbm{1}}
\newcommand{\E}{\mathbb{E}}
\newcommand{\st}{\rm{s.t.}}
\newcommand{\tr}{\mathrm{tr}}
\newcommand{\interior}{\mathrm{int}}

\newtheorem{theorem}{Theorem}[section]
\newtheorem{corollary}[theorem]{Corollary}
\newtheorem{lemma}[theorem]{Lemma}
\newtheorem{proposition}[theorem]{Proposition}
\newtheorem{definition}[theorem]{Definition}
\newtheorem{assumption}[theorem]{Assumption}

\begin{document}
\title{\sf{Piecewise SOS-Convex Moment Optimization and Applications via Exact Semi-Definite Programs  }}

\author{ Q.Y. Huang\footnote{Corresponding author. Department of Applied Mathematics, University of New South Wales, Sydney 2052, Australia. Email: \url{yingkun.huang@unsw.edu.au}. The research of Ms Queenie Yingkun Huang was partially supported by a grant from the Australian Research Council. },  \ V. Jeyakumar\footnote{Department of Applied Mathematics, University of New South Wales, Sydney 2052, Australia. Email: \url{v.jeyakumar@unsw.edu.au}. The research of Prof Vaithilingam Jeyakumar was supported by a grant from the Australian Research Council. } \ and \ G. Li\footnote{Department of Applied Mathematics, University of New South Wales, Sydney 2052, Australia. Email: \url{g.li@unsw.edu.au}. Research of Prof Guoyin Li was supported by a grant from the Australian Research Council }}

\date{\today}
\maketitle

\begin{abstract}
This paper presents exact Semi-Definite Program (SDP) reformulations for infinite-dimensional moment optimization problems involving a new class of piecewise Sum-of-Squares (SOS)-convex functions and projected spectrahedral support sets. These reformulations show that solving a single SDP finds the optimal value and an optimal probability measure of the original moment problem. This is done by establishing an SOS representation for the non-negativity of a piecewise SOS-convex function over a projected spectrahedron. Finally, as an application and a proof-of-concept illustration, the paper also presents numerical results for the Newsvendor and revenue maximization problems with higher-order moments by solving their equivalent SDP reformulations. These reformulations promise a flexible and efficient approach to solving these models. The main novelty of the present work in relation to the recent research lies in finding the solution to moment problems, for the first time, with piecewise SOS-convex functions from their numerically tractable exact SDP reformulations. 
\end{abstract}


\section{Introduction } \label{sec:intro}

Consider the generalized moment problem \begin{align} \label{problem:intro_mmt} \tag{P}
    \min_{\mu \in \mathcal{P}_\Omega}\ \Big\{ \E_\mu^\Omega \left [\min_{k=1,\ldots,r}\max_{\ell=1,\ldots,L} g_\ell^k (\omega) \right ] \; : \; \E_\mu^\Omega [h_j(\omega)] \leq \gamma_j, j=1, \ldots,J \Big \}
\end{align}where $g_\ell^k$, $\ell=1,\ldots,L$, $k=1,\ldots, r$, and $h_j$, $j=1,\ldots,J$, are Sum-of-Squares (SOS)-convex polynomials. The SOS-convexity for polynomials is a numerically tractable relaxation of convexity since checking whether a polynomial is SOS-convex or not can be achieved by solving an SDP \cite{ahmadi2013complete,helton2008structured}. Recent applications of SOS-convexity in optimization may be found in \cite{lasserre2015introduction,jeyakumar2024sum}. The set $\Omega \subset \R^m$ is a convex compact projected spectrahedron and $\mathcal{P}_\Omega$ is the set of probability measures supported on $\Omega$. The expectation of a random variable $\omega$ with respect to the probability measure $\mu \in \mathcal{P}_\Omega$ is denoted by $\E_\mu^\Omega [\omega]$. 

Our model \eqref{problem:intro_mmt} is closely related to generalized moment problems  \cite{lasserre2009moments} where $g_\ell^k$, $\ell=1,\ldots,L$, $k=1,\ldots,r$, and $h_j$, $j=1,\ldots,J$, are arbitrary real-valued continuous functions and the support set $\Omega$ is an arbitrary compact subset of $\R^m$. The dominant technique to reformulate a generalized moment problem as a semi-infinite optimization problem is via the duality theory \cite{shapiro2001duality,zhen2023unified}. However, in general, a semi-infinite optimization problem can be computationally intractable. 

The classical moment problem can be recovered from \eqref{problem:intro_mmt} if $L=1, r=1$, the random variable is univariate, and the constraints are the first $J$ moments \cite{lasserre2009moments}. The multivariate classical moment problem extends to countably many moments, but numerically tractable forms are limited. When $L=1$, $r=1$, $g_1^1$ and $h_j$, $j=1,\ldots,J$, are linear functions and the support $\Omega$ is a spectrahedron, \eqref{problem:intro_mmt} admits an exact SDP reformulation \cite{huang2024distributional}. When $L =1$, $g^{k}_1$, $k=1,\ldots,r$, are linear functions with mean and variance constraints, the support $\Omega$ is an ellipsoid, it has been shown that \eqref{problem:intro_mmt} admits an exact SDP reformulation \cite{delage2010distributionally}. 

When $L=1$, $r=1$, $g_1^1$ and $h_j$, $j=1,\ldots,J$, are any polynomials and the support $\Omega$ is a compact semi-algebraic set, the optimal solution for $\eqref{problem:intro_mmt}$ can be found by solving a convergent hierarchy of SDPs \cite{lasserre2008semidefinite}. Further, polynomial optimization techniques were employed in \cite{klerk2020distributionally} to obtain exact SOS reformulations when $\probset$ contains probability measures whose density functions are SOS polynomials. A review of generalized moment problems and their applications is available in \cite{klerk2019survey}. 

The model problem \eqref{problem:intro_mmt} is of great interest in distributionally robust optimization \cite{zhen2023unified}, where \eqref{problem:intro_mmt} appears as an uncertainty quantification problem to mitigate risks and uncertainties from the input data \cite{zhen2023unified,han2015convex}. For example, in the Newsvendor model of \cite{wiesemann2014distributionally}, the cost of ordering represented as a piecewise linear loss function is minimized to meet the uncertain demand. In the portfolio management problem of \cite{delage2010distributionally}, the profit of the portfolio represented as a piecewise linear utility function is maximized amid random investment returns. 
 
Motivated by the importance of such piecewise functions in optimization and the usefulness of the minimax functions of the form $\min_{k=1,\ldots,r} \max_{\ell=1,\ldots,L} g_\ell^k$ in applications across diverse domains, we introduce a new class of functions, called the piecewise SOS-convex functions. To the best of our knowledge, this notion of piecewise SOS-convex function has so far not been studied in the literature. 

A function $f$ on $\R^n$ is \textit{piecewise SOS-convex} if there exist SOS-convex polynomials $g_\ell^k$, $\ell=1,\ldots,L$, $k=1,\ldots,r$, on $\R^n$ such that $\displaystyle f (v) = \min_{k=1,\ldots,r} \max_{\ell=1,\ldots,L} g_\ell^k (v)$ for all $v \in \R^n$. 

Note that a piecewise SOS-convex function is not necessarily a convex or a differentiable function and it covers a broad range of functions that appear in applications across several domains, as outlined below. Details of their piecewise representations are provided in Appendix \ref{sec:piecewise}.  \begin{itemize}

    \item The truncated $\ell_1$-norm, $f (v) = \min \{1 , \varepsilon |v| \}$, $\varepsilon > 0$, which is non-convex and non-smooth, is used extensively in machine-learning regularization \cite{le2014feature}. 

    \item The piecewise linear function $\displaystyle f(v) = \min_{k=1,\ldots,r} p_k^\top v + q_k$, $p_k \in \R^n, q_k\in \R$, $k=1,\ldots,r$, appears in Newsvendor-type \cite{bertsimas2002relation,goh2010distributionally}, conditional value-at-risk \cite{rockafellar2000optimization}, and portfolio selection \cite{delage2010distributionally} models, and it is a useful approximation for various important utility functions \cite{delage2010distributionally}.

    \item The class of max-SOS-convex functions $\displaystyle f (v) = \max_{\ell=1,\ldots,L} g_\ell (v)$, where $g_\ell$'s are SOS-convex polynomials, is studied in \cite{jeyakumar2014dual}. Further, convex quadratics and convex separable polynomials are SOS-convex \cite{ahmadi2013complete,jeyakumar2012exact}.
    
    \item The class of difference-of-convex functions of the form $\displaystyle f(v) = \max_{\ell=1,\ldots,L} f_\ell (v) - \max_{k=1,\ldots,r} (p_{k}^\top v + q_{k})$, where $f_\ell$'s are SOS-convex polynomials, can be expressed as $\displaystyle f(v) = \min_{k=1,\ldots,r} \max_{\ell=1,\ldots,L} g_\ell^k (v)$ for $g_\ell^k (v) = f_\ell (v) - (p_{k}^\top v + q_{k})$, $\ell=1,\ldots,L$, $k=1,\ldots,r$. These functions have attracted applications in feature selection models of machine-learning \cite{jeyakumar2024sum}. 

    \item The piecewise convex quadratic function \begin{align*}
        f (v) = \begin{cases}
            a (v -b)^2 + c,\ & \text{if}\ 0 \leq v \leq b, \\ c,\ & \text{if}\ v > b. 
        \end{cases}
    \end{align*}with $a,b \geq 0$, $c\in \R$, is used to model offer prices from customers \cite{han2015convex}. 

    \item The Huber loss function with parameter $\varepsilon > 0$, \begin{align*}
        f (v) = \begin{cases}
            \frac{1}{2} v^2,\ & \text{if}\ |v| \leq \varepsilon, \\ \varepsilon |v| - \frac{1}{2} \varepsilon^2,\ & \text{otherwise}, 
        \end{cases}
    \end{align*}is studied in robust statistics \cite{wiesemann2014distributionally}.  
\end{itemize}  

The present study of moment optimization is motivated by two aspects. Firstly, moment problems are numerically challenging due to the presence of infinite-dimensional distributions and multi-dimensional integrals. In response, we examine equivalent Semi-Definite Program (SDP) reformulations of these problems which can be efficiently solved numerically by commonly available software via interior point methods. 

Secondly, the study of piecewise SOS-convexity in \eqref {problem:intro_mmt} was stimulated by the computationally attractive features of SOS-convex polynomial optimization and its applications in many areas. It is known that SOS-convex polynomial optimization problems admit numerically favourable SDP reformulations \cite{jeyakumar2012exact,jeyakumar2015robust,jeyakumar2014dual,lasserre2015introduction}. 

The main goal of this paper is to provide numerically tractable reformulations for the moment problem involving piecewise SOS-convex functions and show how the solution can be recovered from its SDP reformulation. 

\textbf{Our contributions}. Our main contributions are as follows. \begin{enumerate}[label=(\roman*)]

    \item Firstly, we introduce the notion of piecewise SOS-convex functions and establish a new SOS representation of the non-negativity of a piecewise SOS-convex function over a compact projected spectrahedron. This result not only generalizes the known representation result of a non-negative SOS-convex polynomial over a spectrahedron \cite{jeyakumar2012exact} but also provides numerically tractable representations of non-negativity for a broad class of functions that are not necessarily convex. 
    
    \item Secondly, exploiting the representation result, we derive an equivalent numerically tractable SOS optimization reformulation for the generalized moment problem \eqref{problem:intro_mmt}. Moreover, under suitable conditions, we show how to recover an optimal probability measure for \eqref{problem:intro_mmt} by solving a single SDP. This extends the corresponding result of Lasserre \cite{lasserre2009convexity} on SOS-convex polynomial optimization.
      
    \item Finally, as an application and a proof-of-concept illustration, we present numerical results for two important practical models, i.e., the Newsvendor problem \cite{guo2022unified} with higher-order moments and the revenue maximization problem \cite{han2015convex} with a higher-order utility function, by solving their equivalent SDP reformulations.    
\end{enumerate}

The novelty of the present work in relation to recent research in SOS-convexity and moment optimization is the derivation of numerically tractable exact SDP reformulations for moment problems involving, {\it for the first time, piecewise SOS-convex functions}. The main innovation is the combined utilization of powerful tools from real algebraic geometry, convex analysis, and SOS polynomials to produce the exact reformulations. We employ the key computationally favourable features of SOS-convexity within an infinite-dimensional setting and exploit the geometry of the projected spectrahedra to facilitate the SDP reformulations. 

Our approach makes use of the piecewise structure of the functions for transforming infinite-dimensional problems into finite-dimensional numerically tractable problems. Moreover, the projected spectrahedron covers a broad class of semi-algebraic sets \cite{lasserre2015introduction}, such as the spectrahedra, ellipsoids, and boxes, used in robust optimization \cite{ben2009robust}. 

The paper is organized as follows. Section \ref{sec:sos} provides SOS representation results for piecewise SOS-convex functions. Section \ref{sec:mmt} presents the main theorems for the SOS reformulation of \eqref{problem:intro_mmt} and the optimal solution recovery of \eqref{problem:intro_mmt}. Section \ref{sec:applications} describes applications to the Newsvendor and revenue maximization models. Section \ref{sec:conclusion} concludes with discussions on potential future work. The appendices provide technical details related to explicit piecewise representations for some piecewise SOS-convex functions (Appendix \ref{sec:piecewise}), infinite-dimensional conic duality (Appendix \ref{sec:conic}), and SDP representations of SOS optimization problems (Appendix \ref{sec:sos-sdp}). 

\section{Non-Negativity of Piecewise SOS-Convex Functions } \label{sec:sos}

In this section, we outline SOS representations for a class of piecewise functions involving SOS-convex polynomials. This will play a vital role in reformulating the generalized moment problem \eqref{problem:intro_mmt} as a numerically tractable SOS optimization problem. 

We begin with some preliminaries on polynomials. Denote $\R^m$ the Euclidean space of dimension $m$ and $\R^m_+$ the non-negative orthant of $\R^m$. The standard inner product is $a^\top b$ for $a,b \in \R^m$. Let $\R [v]$ be the space of polynomials with real coefficients over $v \in \R^m$. A polynomial $f \in \R[v]$ is called a Sum-of-Squares (SOS) polynomial if there exist polynomials $f_j \in \R[v], j = 1,\ldots,J$, such that $f = \sum_{j=1}^J f_j^2$. For a polynomial $f \in \R[v]$, we use ${\rm deg}\, f$ to denote its degree. We also use $\Sigma_d^2 (v)$ to denote the set of all SOS polynomials $f$ of degree at most $d$ with respect to the variable $v\in \R^m$. Next, we recall the definition of SOS-convex polynomials.

\begin{definition}[\textbf{SOS-convex polynomial} \cite{ahmadi2013complete,helton2010semidefinite}]

    A polynomial $f \in \R[v]$ is SOS-convex if its Hessian $H(v)$ is an SOS matrix polynomial, that is, if there exists an $(s\times m)$ polynomial matrix $P(v)$ for some $s \in \mathbb N$ such that $H(v) = P (v)^\top P(v)$. 
\end{definition}

Several equivalent conditions for SOS-convexity are available in \cite{ahmadi2013complete}. For instance, 
$f \in \R[v]$ is SOS-convex whenever the polynomial $g(v_1, v_2) = f(v_1) - f(v_2) - \nabla f(v_2)^\top (v_1 - v_2)$ on $\R^m \times \R^m$ is an SOS polynomial with respect to the variable $(v_1, v_2)$. 

Any SOS-convex polynomials are convex polynomials, but the converse is not true \cite{ahmadi2013complete}. In other words, the class of SOS-convex polynomials is a proper subclass of convex polynomials. Moreover, the class of SOS-convex polynomials covers affine functions, convex quadratics, and convex separable polynomials, while SOS-convex polynomials can be non-quadratic or non-separable in general \cite{jeyakumar2012exact}. 

Now, we formally define the notion of a piecewise SOS-convex function. 

\begin{definition}[{\bf Piecewise SOS-convex function}] \label{defn:piecewise-sos}
    We call a function $f$ on $\R^m$ piecewise SOS-convex if there exist SOS-convex polynomials $g_\ell^k$, $\ell=1,\ldots,L$, $k=1,\ldots,r$, on $\R^m$ such that $\displaystyle f (v) = \min_{k=1,\ldots,r} \max_{\ell=1,\ldots,L} g_\ell^k (v)$ for all $v \in \R^m$. 
\end{definition}

The following proposition links non-negative SOS-convex polynomials to SOS polynomials which will be deployed in the tractable representation of piecewise SOS-convex functions (Theorem \ref{thm:max_sos_conv_spect}). 

\begin{proposition} (see \cite{helton2010semidefinite} and \cite[Corollary~2.1]{jeyakumar2014dual}). \label{prop:nn_sos_conv}
    Let $f \in \R[v]$ be a non-negative SOS-convex polynomial. Then, $f$ is an SOS polynomial. 
\end{proposition}

Denote $e_i \in \R^m$ the $i$-th standard basis vector of $\R^m$, $i=1,\ldots,m$. The simplex in $\R^L$ is defined as $\Delta = \{\delta \in \R_+^{L} : \sum_{\ell=1}^{L} \delta_{\ell} = 1\}$. Let $\mathbb{S}^\nu$ be the set of $(\nu \times \nu)$ symmetric matrices. The trace product on $\mathbb S^\nu$ is $\tr(MN)$ for any $M,N \in \mathbb S^\nu$. A matrix $M \in \mathbb S^\nu$ is positive semi-definite, denoted as $M \succeq 0$ (resp. positive definite, denoted as $M \succ 0$), if $x^\top M x \geq 0$  for all $x\in \R^\nu$ (resp. $x^\top M x >0$ for all $x \in \R^\nu$, $x \neq 0$). Let $\mathbb{S}_+^\nu$ be the cone of symmetric $(\nu \times \nu)$ positive semi-definite matrices. 

A projected spectrahedron \cite{helton2008structured}, or a spectrahedral shadow \cite{netzergeometry}, takes the form \begin{align} \label{eqn:lmi} \tag{LMI}
    \Omega = \left \{v \in \R^m : \exists \xi \in \R^N, F_0 + \sum_{i=1}^m v_i F_i + \sum_{t=1}^N \xi_t M_t \succeq 0 \right \},  
\end{align}for some $F_i \in \mathbb S^\nu, i=0,1,\ldots,m$, and $M_t \in \mathbb S^\nu, t = 1,\ldots,N$. This set will later be used as a support set for probability measures. This set is versatile as any semi-algebraic set $\{v\in \R^m : f_i (v) \leq 0, i=1,\ldots,n \}$ described by SOS-convex polynomials $f_i, i =1,\ldots,n$, can be represented as a projected spectrahedron \cite{helton2008structured}. Moreover, spectrahedra, ellipsoids, and boxes are special forms of projected spectrahedra. 

The following representation of a piecewise SOS-convex function extends the result in \cite{jeyakumar2012exact} for an SOS-convex polynomial over a spectrahedron and in \cite{jeyakumar2014dual} for a max-SOS-convex polynomial over a subclass of projected spectrahedron. 

\begin{theorem}[\textbf{SOS representation for non-negative piecewise SOS-convex functions}] \label{thm:max_sos_conv_spect} 
    Let $\Omega \subset \R^m$ be given as in \eqref{eqn:lmi}, $g_\ell^k, \ell = 1,\ldots,L, k=1,\ldots,r$, be SOS-convex polynomials on $\R^m$. Assume that $\Omega$ is compact and there exist $\bar  v\in \R^m$ and $\bar  \xi\in \R^N$ such that $\displaystyle F_0 + \sum_{i=1}^m \bar  v_i F_i + \sum_{t=1}^N \bar  \xi_t M_t \succ 0$. Then, the following statements are equivalent. \begin{enumerate}
        \item $v\in \Omega \implies \displaystyle \min_{k=1,\ldots,r} \max_{\ell=1,\ldots,L} g_\ell^k (v)  \geq 0$. 

        \item For each $k=1,\ldots,r$, there exist $Z_k \in \mathbb S_+^\nu$ and $\delta^k  = (\delta_1^k, \ldots, \delta_L^k) \in \Delta$ such that \begin{align*}
            \begin{cases}
                \displaystyle \sum_{\ell=1}^L \delta_\ell^k g_\ell^k (v)  - \tr (Z_k F_0) - \sum_{i=1}^m v_i \tr (Z_k F_i) \in \Sigma_d^2 (v),\ k=1,\ldots,r, \\ \tr(Z_k M_t) = 0,\ t = 1,\ldots,N, \; k=1,\ldots, r.
            \end{cases}            
        \end{align*}where $d$ is the smallest even integer larger than $\displaystyle  \max_{k=1,\ldots,r} \max _{\ell = 1,\ldots,L} \deg g_\ell^k$. 
    \end{enumerate}    
\end{theorem}

\begin{proof}
    
    [\textit{(ii)} $\implies$ \textit{(i)}]. Assume that there exist $Z_k \in \mathbb S_+^\nu$ and $\delta^k  \in \Delta$, $k=1,\ldots,r$, such that \begin{align*} \begin{cases}
        \displaystyle \sum_{\ell=1}^L   \delta_\ell^k g_\ell^k (v) - \tr(Z_k F_0) - \sum_{i=1}^m v_i \tr (Z_k F_i) \in \Sigma_d^2 (v),\ k=1,\ldots,r, \\ \tr(Z_k M_t) = 0,\ t = 1,\ldots,N,\ k=1,\ldots,r.
    \end{cases} \end{align*} 
    
    Fix $k\in \{1,\ldots,r\}$. Since an SOS polynomial always takes non-negative values, it follows that \begin{align*}
        \sum_{\ell=1}^L \delta_\ell^k g_\ell^k (v) - \tr(Z_k F_0) - \sum_{i=1}^m v_i \tr (Z_k F_i) \geq 0,\ \text{for all}\ v \in \R^m.
    \end{align*}
    
    Further, $\tr (Z_k F_0) + \sum_{i=1}^m v_i \tr (Z_k F_i) = \tr ( Z_k (F_0 + \sum_{i=1}^m v_i F_i + \sum_{t=1}^N \xi_t M_t) )$ for any $\xi \in \R^N$. This implies that, for all $(v,\xi) \in \R^m \times \R^N$, \begin{eqnarray}\label{eq:use1}
       \sum_{\ell=1}^L   \delta_\ell^k g_\ell^k (v)  \geq \tr\Big (Z_k \Big (F_0+ \sum_{i=1}^m v_i F_i + \sum_{t=1}^N \xi_t M_t \Big ) \Big ). 
    \end{eqnarray}

    Since $(\delta_1^k,\ldots,\delta_L^k) \in \Delta$, we have $\displaystyle \max_{\ell=1,\ldots,L} g_\ell^k (v) \geq \sum_{\ell=1}^L \delta_\ell^k g_\ell^k (v)$. Further, for an arbitrary ${v} \in \Omega$, there exists $ \xi \in \R^N$ such that $\displaystyle F_0 + \sum_{i=1}^m v_i F_i + \sum_{t=1}^N \xi_t M_t \in \mathbb S_+^\nu$. Also, since $Z_k\in \mathbb S_+^\nu$, we have $\displaystyle \tr \Big (Z_k \Big (F_0  +\sum_{i=1}^m v_i F_i + \sum_{t=1}^N \xi_t M_t \Big ) \Big ) \geq 0$. It follows from Eqn. \eqref{eq:use1} that $\displaystyle \max_{\ell=1,\ldots,L} g_\ell^k (v) \geq0$ for any $v\in \Omega$. This holds for any $k\in \{ 1,\ldots,r \}$, so $\displaystyle \min_{k=1, \ldots,r } \max_{\ell=1,\ldots,L} g_\ell^k (v) \geq 0$ for any $v \in \Omega$. 

    [\textit{(i)} $\implies$ \textit{(ii)}]. Assume Statement \textit{(i)} is true. Fix $k \in \{1,\ldots,r\}$. We have \begin{align} \label{eqn:max_sos_spect_pf} 
        \min_{v\in \Omega} \max_{\ell=1,\ldots,L} g_\ell^k (v) = \min_{v\in \Omega} \max_{\delta\in \Delta} \sum_{\ell=1}^L \delta_\ell g_\ell^k (v) \geq 0.
    \end{align}

    Define a bifunction $P_k : \Omega \times \Delta \to \R$ by $P_k (v,\delta) = \sum_{\ell=1}^L \delta_{\ell} g_\ell^k (v)$. Direct verification shows that $P_k (v, \cdot)$ is linear for each fixed $v\in \Omega$ and $P_k (\cdot, \delta)$ is convex continuous for each fixed $\delta \in \Delta$. Note that both $\Omega \subseteq \R^m$ and $\Delta \subseteq \R^L$ are convex compact sets. The Convex-Concave Minimax Theorem \cite[cf.][Theorem~2.10.2]{zalinescu2002convex} gives us that \begin{align*}
        \min_{v\in \Omega} \max_{\delta \in \Delta} P_k (v,\delta) = \max_{\delta \in \Delta} \min_{v\in \Omega} P_k (v,\delta).
    \end{align*}
    
    Thus, there exists $\delta^k \in \Delta$ such that $\displaystyle \max_{\delta \in \Delta} \min _{v\in \Omega} P_k (v,\delta) = \min_{v\in \Omega} P_k (v,\delta^k ) = \min_{v\in \Omega} \sum_{\ell=1}^L \delta_\ell^k g_\ell^k (v)$.
        
    Observe that 
    \begin{align}
        \min_{v\in \Omega} \sum_{\ell=1}^L \delta_\ell^k  g_\ell^k (v) = & \min_{v\in \R^m} \left \{ \sum_{\ell=1}^L \delta_\ell^k  g_\ell^k (v)   :    \, F_0 + \sum_{i=1}^m v_i F_i + \sum_{t=1}^N \xi_t M_t \succeq 0  \mbox{ for some } \xi \in \R^N \right\} \nonumber \\ = & \min_{(v,\xi) \in \R^m\times \R^N}\left \{ \sum_{\ell=1}^L \delta_\ell^k g_\ell^k (v) : F_0 + \sum_{i=1}^m v_i F_i + \sum_{t=1}^N \xi_t M_t \succeq 0 \right \}. \label{eqn:maxl_maxdelta}
    \end{align}
    
    By the Lagrangian Duality under the strict feasibility assumption on $\Omega$ \cite[Theorem~2.9.2]{zalinescu2002convex}, there exists $Z_k \in\mathbb S_+^\nu$ such that \begin{align*}
        \min \eqref{eqn:maxl_maxdelta} = \min_{\substack{ (v, \xi) \in \R^m \times \R^N}} \left \{ \sum_{\ell=1}^L \delta_\ell^k  g_\ell^k (v)  - \tr(Z_k F_0) - \sum_{i=1}^m v_i \tr(Z_k F_i) - \sum_{t=1}^N \xi_t \tr(Z_k M_t) \right \}.
    \end{align*}

    Hence, Eqn. \eqref{eqn:max_sos_spect_pf} implies that there exist $Z_k \in \mathbb S_+^\nu$ and $\delta^k \in \Delta$ such that \begin{align} \label{eqn:nonneg_sos}
        \sum_{\ell=1}^L \delta_\ell^k g_\ell^k (v) - \tr(Z_k F_0) - \sum_{i=1}^m v_i \tr(Z_k F_i) - \sum_{t=1}^N \xi_t \tr(Z_k M_t) \geq 0,
    \end{align}for all $(v,\xi) \in \R^m \times \R^N$. Letting $\xi=0 \in \R^N$ in Eqn. \eqref{eqn:nonneg_sos} gives \begin{align*}
        \sigma_k (v) := \sum_{\ell=1}^L \delta_\ell^k g_\ell^k (v) - \tr(Z_k F_0) - \sum_{i=1}^m v_i \tr(Z_k F_i)  \geq 0, \; \mbox{for all} \; v\in \R^m.
    \end{align*} 
    
    Since $\sigma_k$ is a sum of finitely many SOS-convex polynomials, $\sigma_k$ is a non-negative SOS-convex polynomial. By Proposition \ref{prop:nn_sos_conv}, $\sigma_k$ is thus an SOS polynomial. Moreover, by fixing $v=\widehat {v}$ for some $\widehat {v} \in \R^m$, Eqn. \eqref{eqn:nonneg_sos} further implies that $\sigma_k (\widehat  v) - \sum_{t=1}^N \xi_t \tr(Z_k M_t) \geq 0$ for any $\xi\in \R^N$ where $\sigma_k(\widehat {v}) \geq 0$ is fixed. This forces $\tr(Z_k M_t) = 0$ for all $t = 1,\ldots,N$. The conclusion now follows because $k$ is chosen arbitrarily from $\{1,\ldots,r\}$.
\end{proof}

\section{Piecewise SOS-Convex Moment Problems } \label{sec:mmt}

In this section, we focus on reformulating {\it a generalized moment optimization problem with piecewise SOS-convex functions} as an SDP. 

We begin by fixing some terminologies. Let $Y$ and $Y'$ be real topological vector spaces. They are paired if a bilinear form $\langle \cdot, \cdot \rangle : Y' \times Y \to \R$ is defined. The cone generated by any set $S\subseteq Y$ and the interior of $S$ are denoted as cone$(S)$ and $\interior(S)$ respectively. The (positive) polar of a cone $S \subseteq Y$ is defined as $S^+ = \{y' \in Y' : \langle y', y \rangle \geq 0, \forall  y \in S\}$. See Appendix \ref{sec:conic} for details on strong conic duality. 

For a non-empty compact set $\Omega \subset \R^m$, $\mathcal{C}_\Omega$ denotes the vector space of all continuous real-valued functions on $\Omega$ equipped with the supremum norm. Denote $\mathcal{B}$ the Borel $\sigma$-algebra of $\Omega$. {\it From this section onwards, we let $X$ be the vector space of (finite signed regular) Borel measures on $(\Omega, \mathcal{B})$}. It is known that $X$ can be identified as the topological dual space of $\mathcal{C}_\Omega$, i.e., $X = \mathcal{C}_\Omega^*$. We equip $X$ with the weak$^*$ topology. For any point $v\in \Omega$, the Dirac measure $\1_v$ which takes a mass of $1$ at the point $v\in \Omega$ and $0$ otherwise belongs to $X$ \cite{shapiro2001duality}. 

Consider the following generalized moment problem with piecewise SOS-convex functions:  \begin{align} \label{problem:uq} \tag{P} 
    \min_{\mu \in X}\ & \E_\mu^\Omega \Big [ \min_{k=1,\ldots,r} \max_{\ell=1,\ldots,L} g_\ell^k (\omega) \Big ] \\ \st\ & \E_\mu^\Omega \left [h_j (\omega) \right ] \leq \gamma_j,\ j=1,\ldots, J, \nonumber \\ & \E_\mu^\Omega [1] = 1,\ \mu \succeq_\B 0, \nonumber
\end{align}where $\gamma_j \in \R$, $j=1,\ldots,J$, and $h_j$, $j = 1,\ldots,J$, and $g_\ell^k$, $\ell =1,\ldots, L$, $k=1,\ldots,r$, are SOS-convex polynomials on $\R^m$. Define $g: \R^m \rightarrow \R$, $\displaystyle g(v)= \min_{k=1,\ldots,r} \max_{\ell=1,\ldots,L} g_\ell^k (v)$ the piecewise SOS-convex function, and note that $g \in \mathcal{C}_\Omega$ and $h_j \in \mathcal{C}_\Omega$, $j=1,\ldots,J$. We let $X'$ be the vector space generated by $g$, $h_j, j=1,\ldots,J$, and the constant function of $1$. Note that $X' \subseteq \mathcal{C}_\Omega$, and $X'$ is equipped with the supremum norm. Following \cite{shapiro2001duality}, $X$ and $X'$ are paired spaces, and the continuous bilinear form between them is given by \begin{align*}
    \langle  f,\mu  \rangle = \E_\mu^\Omega [ f(\omega)]:= \int_\Omega f (\omega)\ d \mu (\omega), \mbox{ for all } f\in X' \mbox{ and } \mu\in X,
\end{align*}where $\E_\mu^\Omega [f (\omega)]$ refers to the expectation of the random variable $f(\omega)$ with respect to the measure $\mu$ supported on $\Omega$. Denote the set of probability measures by $\mathcal{P}_\Omega = \{ \mu \in X : \E_\mu^\Omega [1] = 1, \mu \succeq_\B 0 \}$, where $\mu\succeq_\B 0$ means that $\mu (A) = \mu (\omega \in A) \geq 0$ for all $\mathcal{B}$-measurable sets $A$. 

The moment problem \eqref{problem:uq} is known as an uncertainty quantification problem in the distributionally robust optimization literature \cite{hanasusanto2015distributionally}, and the feasible set $\{ \mu\in \mathcal{P}_\Omega : \E_\mu^\Omega [h_j (\omega)] \leq \gamma_j, j=1,\ldots,J\}$ is called a moment ambiguity set. 

We assume that the support set $\Omega$ is a \textit{(convex) projected spectrahedron}: \begin{align} \tag{LMI}
    \Omega = \Big \{v \in \R^m : \exists \xi \in \R^N,\ F_0 + \sum_{i=1}^m v_i F_i + \sum_{t=1}^N \xi_t M_t \succeq 0 \Big \},
\end{align} for some $F_i \in \mathbb S^\nu, i=0,1,\ldots,m$, and $M_t \in \mathbb S^\nu, t = 1,\ldots,N$. 
We associate \eqref{problem:uq} with the following SOS optimization problem \begin{align} \label{problem:uq_dual}    \tag{D} 
   \max_{\substack{ \lambda \in \R^{J}_+\times \R \\ Z_k \in \mathbb S_+^\nu, \delta^k  \in \Delta \\ k=1,\ldots, r}} \ & -\sum_{j=1}^{J} \lambda_{j} \gamma_j - \lambda_{J+1} \\ \st\hspace{10mm} & \sum_{\ell=1}^{L} \delta_\ell^k  g_\ell^k (v) + \sum_{j=1}^{J} \lambda_{j} h_j (v) + \lambda_{J+1}  - \tr(Z_kF_0) -  \sum_{i=1}^m v_i \tr(Z_kF_i) \in \Sigma_d^2 (v),\ k=1,\ldots,r, \nonumber  \\ & \tr(Z_k M_t) = 0, \ t = 1,\ldots,N, \  k=1,\ldots,r, \nonumber 
\end{align}where $\delta^k  = (\delta_1^k, \ldots, \delta_L^k) \in \Delta = \{\delta \in \R^L_+ : \sum_{\ell=1}^L \delta_\ell = 1\}$, $\Sigma_d^2 (v)$ is the set of SOS polynomials of degree at most $d$ with respect to the variable $v \in\R^m$, and $d$ is the smallest even integer with \begin{align*}
    d\geq \max \left \{\max_{k=1,\ldots,r} \max_{\ell = 1,\ldots,L} \deg g_\ell^k, \max_{j=1,\ldots,J}  \deg h_j \right \}.
\end{align*}

The SOS program \eqref{problem:uq_dual} can be equivalently rewritten as a Semi-Definite Program (SDP), and thus can be efficiently solved via existing SDP software. We refer the readers to Appendix \ref{sec:sos-sdp} for the procedure. 

The following useful characterization states that the polar cone of the cone of non-negative measures is the convex cone of functions in $X'$ that are non-negative on $\Omega$. This characterization is known (e.g., \cite[Examples~2.37, 2.122]{bonnans2013perturbation}, \cite[Lemma~3.1]{huang2024distributional}), but we provide the proof here for the sake of self-containment. 

Note that, for the convex cone of non-negative measures $C = \{\mu \in X : \mu \succeq_\B 0 \}$, we denote $C^+ = \{f \in X' : \langle  f, \mu \rangle \geq 0, \forall \mu \in C \}$. Here,  $C^+$ is the intersection of the topological polar cones of $C$ with the subspace $X'$ of the vector space of continuous real-valued functions $\mathcal{C}_\Omega$.

\begin{lemma}[{\bf Characterization of cone of non-negative measures}]\label{lemma:max_supp}
    Let $\Omega \neq \emptyset$ be any convex compact subset of $\R^m$, $C = \{\mu \in X: \mu \succeq_\B 0\}$, and $f\in X'$. Then, $f \in C^+$ if and only if $\min_{v \in \Omega} f(v) \geq 0$.
\end{lemma}

\begin{proof}
    Suppose that $f \in C^+$. Then, $\langle  f, \mu \rangle \geq 0$ for all $\mu \in C$. Note that the Dirac measure $\1_v$ at any point $v \in \Omega$ belongs to $C$. This gives $f(v) = \langle  f, \1_v \rangle \geq 0$ for any $v \in \Omega$. Hence, $\displaystyle \min_{v \in \Omega} f(v) \geq 0$. 

    Conversely, suppose that $\displaystyle \min_{v\in \Omega} f(v) \geq 0$. Now, for any $\mu \in \mathcal{P}_\Omega = \{ \mu \in X : \int_\Omega 1 d\mu (\omega) = 1, \mu \succeq_\B 0 \}$, we have \begin{align*}
        \int_\Omega f (\omega)\ d\mu (\omega) \geq \int_\Omega \left (\min_{v\in \Omega} f (v) \right ) d\mu (\omega) = \left (\min_{v \in \Omega} f (v) \right ) \int_\Omega 1 d\mu (\omega) = \min_{v \in \Omega} f (v).
    \end{align*}
    
    This implies $\langle  f,\mu \rangle \geq 0$ for any $\mu \in \mathcal{P}_\Omega$. Since $C =$ cone$(\mathcal{P}_\Omega)$, we see $f \in C^+$. 
\end{proof}

Now, we show that \eqref{problem:uq_dual} is an exact SOS reformulation for \eqref{problem:uq} in the sense that $\min \eqref{problem:uq} = \max \eqref{problem:uq_dual}$.

\begin{assumption}[{\bf Interior point condition for moment problem}] \label{assump:interior}
    \begin{align*}
    (\gamma_1, \ldots, \gamma_J, 1) \in \interior \bigg\{\{ ( \langle  h_1, \mu \rangle, \ldots, \langle  h_J, \mu \rangle, \langle  1, \mu \rangle) : \mu \in X, \mu \succeq_\B0 \} + \R_+^{J} \times \{0\} \bigg\}.
\end{align*} 
\end{assumption}

This interior point condition is commonly used in the moment problem literature and can be satisfied in many situations. See \cite{shapiro2001duality,xu2018distributionally} for discussions. 

\begin{assumption} \label{assump:uq}
    The projected spectrahedron $\Omega$ as in \eqref{eqn:lmi} is compact and there exist $\bar  v \in \R^m$ and $\bar \xi \in \R^N$ such that $F_0 + \sum_{i=1}^m \bar  v_i F_i + \sum_{t=1}^N \bar  \xi_t M_t\succ 0$. For the problem \eqref{problem:uq}, the polynomials $g_\ell^k$, $\ell=1,\ldots,L$, $k=1,\ldots,r$, and $h_j$, $j=1,\ldots,J$, are SOS-convex polynomials on $\R^m$. 
\end{assumption}

\begin{theorem}[\textbf{Exact SOS program for piecewise SOS-convex moment optimization}] \label{thm:uqd_max_lin}
    Assume that \eqref{problem:uq} admits a minimizer and Assumptions \ref{assump:interior}-\ref{assump:uq} are satisfied. Then, $\min \eqref{problem:uq} = \max \eqref{problem:uq_dual}$. 
\end{theorem}

\begin{proof} 

     \textbf{[Conic duality]}. Let $g(v) = \displaystyle  \min_{k=1,\ldots,r} \max_{\ell=1,\ldots,L} g_\ell^k (v)$, $b =(\gamma_1,\ldots,\gamma_{J}, 1)$, $K = \R_+^{J} \times \{0\}$, and $C = \{\mu \in X : \mu \succeq_\B 0\}$. Note that $g \in X'$, and $C$ and $K$ are convex cones that are closed in the respective topologies. Define a continuous linear map $A : X \to \R^{J+1}$ by $A \mu=  (\langle h_1, \mu \rangle, \ldots, \langle h_J, \mu \rangle, \langle  1, \mu \rangle)$. Then, \begin{align*}
         \min \eqref{problem:uq} = \min_{\mu \in C}\{\langle g, \mu \rangle : -A \mu +b  \in K \}.
     \end{align*}
     
     The adjoint mapping $A^*: \R^{J+1} \mapsto X'$ is given by $ A^* \lambda = \sum_{j=1}^{J} \lambda_j h_j + \lambda_{J+1} \cdot 1$, and the (positive) polar cone of $K$ is $K^+ = \R_+^{J} \times \R$. Recall that $C^+ = \{f\in X' : \langle f,\mu \rangle \geq 0, \forall \mu \in C\}$. Then, Assumption \ref{assump:interior} can be equivalently rewritten as $b \in \interior (A (C) + K)$. It follows from strong conic linear duality (Corollary \ref{coro:int_pt}) that \begin{eqnarray} \label{eqn:conic-uqdual}
        \min \eqref{problem:uq} = \max \left \{ - \sum_{j=1}^{J} \lambda_j \gamma_j - \lambda_{J+1} : g + \sum_{j=1}^{J} \lambda_j h_j + \lambda_{J+1} \cdot 1 \in C^+, \lambda\in K^+ \right \},
    \end{eqnarray}and the maximum of the problem in Eqn. \eqref{eqn:conic-uqdual} is attained. 

    \textbf{[Equivalent SOS representation of conic constraint]}. This means, there exists $\lambda \in \R^{J}_+ \times \R$ with $\displaystyle g+ \sum_{j=1}^{J} \lambda_j h_j + \lambda_{J+1} \cdot 1 \in C^+$ such that $\displaystyle \min\eqref{problem:uq} = -\sum_{j=1}^{J} \lambda_j \gamma_j - \lambda_{J+1}$. By Lemma \ref{lemma:max_supp}, $g + \sum_{j=1}^{J} \lambda_j h_j + \lambda_{J+1} \cdot 1 \in C^+$ if and only if \begin{align} \label{eqn:sos_rep}
        \min_{v \in \Omega} \left \{ \min_{k=1,\ldots,r} \max_{\ell=1,\ldots,L} g_\ell^k (v) + \sum_{j=1}^{J} \lambda_j h_j(v) + \lambda_{J+1} \right \} \geq 0. 
    \end{align}
    
    The polynomial $g_\ell^k +\sum_{j=1}^J \lambda_j h_j + \lambda_{J+1}$ is SOS-convex for each $\ell = 1,\ldots, L$, $k=1,\ldots,r$. By Theorem \ref{thm:max_sos_conv_spect}, Eqn. \eqref{eqn:sos_rep} is equivalent to, for each $k=1,\ldots,r$, \begin{align*} 
        \begin{cases} \displaystyle
            \sum_{\ell=1}^L \delta_\ell^k g_\ell^k (v) + \sum_{j=1}^{J} \lambda_j h_j(v) + \lambda_{J+1} -\tr(Z_k F_0) - \sum_{i=1}^m v_i \tr(Z_k F_i)\in \Sigma_d^2 (v), \\ \tr(Z_k M_t) = 0, \ t = 1,\ldots,N, 
        \end{cases} 
    \end{align*}for some $Z_k \in \mathbb S_+^\nu$ and $\delta^k  \in \Delta$. Therefore, $\min \eqref{problem:uq} \leq \max \eqref{problem:uq_dual}$. 
    
    \textbf{[Weak duality]}. Take any feasible point $\lambda \in \R^{J}_+ \times \R$, $Z_k \in \mathbb S_+^\nu$, $\delta^k  \in \Delta, k=1,\ldots, r$, for \eqref{problem:uq_dual}. Theorem \ref{thm:max_sos_conv_spect} gives $\displaystyle \min_{k=1,\ldots,r} \max_{\ell=1,\ldots, L} \Big \{ g_\ell^k (v) + \sum_{j=1}^J \lambda_j h_j(v) + \lambda_{J+1} \Big \}  \geq0$ for all $v \in \Omega$, and thus, \begin{align*} 
    \min_{k=1,\ldots,r} \max_{\ell=1,\ldots, L} g_\ell^k (v) \geq - \sum_{j=1}^J \lambda_j h_j(v) - \lambda_{J+1}, \ \text{for all}\ v \in \Omega.
    \end{align*}

    Take any feasible point $\mu$ for \eqref{problem:uq}, then, $\displaystyle \E_\mu ^\Omega \Big [ \min_{k=1,\ldots,r} \max_{\ell=1,\ldots, L} g_\ell^k (\omega) \Big ] \geq \E_\mu^\Omega \Big  [ - \sum_{j=1}^J \lambda_j h_j(\omega) - \lambda_{J+1} \Big ]$. Further, $\displaystyle -\sum_{j=1}^J \lambda_j \gamma_j - \lambda_{J+1} \leq \E_\mu^\Omega \Big [-\sum_{j=1}^J \lambda_j h_j(\omega) - \lambda_{J+1} \Big ]$. Thus, $\displaystyle \E_\mu ^\Omega \Big [ \min_{k=1,\ldots,r} \max_{\ell=1,\ldots, L} g_\ell^k (\omega) \Big ] \geq -\sum_{j=1}^J \lambda_j \gamma_j - \lambda_{J+1}$, implying $\min \eqref{problem:uq} \geq \max\eqref{problem:uq_dual}$. Together, $\min \eqref{problem:uq} = \max \eqref{problem:uq_dual}$.
\end{proof}

The SOS optimization problem \eqref{problem:uq_dual} can be expressed equivalently as the following SDP (Proposition \ref{prop:sos-sdp}), \begin{align} \label{problem:uq_sdp} \tag{R}
    \max_{\substack{ \lambda \in \R^J_+ \times \R,\delta^k \in \Delta \\ Z_k \in \mathbb S_+^\nu, Q_k \in \mathbb S_+^{s (m, d/2) }  \\ k=1,\ldots,r } }\ & -\sum_{j=1}^J \lambda_j \gamma_j - \lambda_{J+1} \\ \st \hspace{7mm} \ & \sum_{\ell=1}^L \delta_\ell^k (g_\ell^k)_\alpha + \sum_{j=1}^J \lambda_j (h_j)_\alpha + \lambda_{J+1} (1)_\alpha - \tr(Z_k F_\alpha) = \tr(Q_k B_\alpha), \nonumber \\  & \hspace{20mm} \text{for all}\ \alpha \in \mathcal{N},\ k=1,\ldots,r, \nonumber \\ & \tr(Z_k M_t) = 0,\ t=1,\ldots,N,\ k=1,\ldots,r, \nonumber
\end{align}where $s(m,d) = \binom{m+d}{d}$, $\mathcal{N} = \{\alpha = (\alpha_1, \ldots, \alpha_m) \in (\{ 0\} \cup \mathbb N)^m : \sum_{i=1}^m \alpha_i \leq d \}$, and $B_\alpha \in \mathbb S^{s(m, d/2)}$, $\alpha \in \mathcal{N}$, are matrices given in Appendix \ref{sec:sos-sdp}. The notation $(f)_\alpha$ of any polynomial $f\in \R[v]$ of degree at most $d$ refers to the $\alpha$-th coefficient of $f$. In particular, $(1)_\alpha$ refers to the $\alpha$-th coefficient of the constant function $1$, that is, $(1)_\alpha = 1$ when ${\alpha = 0}$, and $(1)_\alpha = 0$ otherwise. The matrices are $F_0 := F_0 \in \mathbb{S}^{\nu}$, $F_{e_i} := F_i \in \mathbb{S}^{\nu}$, $i=1,\ldots,m$, and $F_\alpha$, $\alpha \in \mathcal{N} \setminus \{0, e_1,\ldots,e_m\}$, is the zero matrix. 

The dual problem of \eqref{problem:uq_sdp}, which is also an SDP, is given by  \begin{align} \label{problem:uq_sdp_recovery} \tag{S}
    \min_{\substack{ y^k \in \R^{s(m, d)}, \xi^k \in \R^N\\ z^k \in \R, k=1,\ldots,r } }\ & \sum_{k=1}^r z^k \\ \st \hspace{13mm} & \sum_{k=1}^r \sum_{\alpha \in \mathcal{N}} y_\alpha^k (h_j)_\alpha \leq \gamma_j,\ j=1,\ldots,J, \label{eqn:recovery1} \\ & \sum_{\alpha \in \mathcal{N}} y_\alpha^k (g_\ell^k)_\alpha \leq z^k,\ \ell =1,\ldots,L,\ k=1,\ldots,r, \label{eqn:recovery2} \\ & y_0^k F_0 + \sum_{i=1}^m y_{e_i}^k F_i + \sum_{t=1}^N \xi_t^k M_t \succeq 0,\ k=1,\ldots, r, \label{eqn:recovery3} \\ & \sum_{k=1}^r y_0^k = 1,\ \sum_{\alpha \in \mathcal{N}} y_\alpha^k B_\alpha \succeq 0,\ k=1,\ldots,r, \label{eqn:recovery4}
\end{align}

Denote the minimum value of the problem in \eqref{problem:uq_sdp_recovery} as $\min \eqref{problem:uq_sdp_recovery}$. Sufficient conditions for the strong duality between the conic programs \eqref{problem:uq_sdp} and \eqref{problem:uq_sdp_recovery} are standard \cite[cf.][]{bonnans2013perturbation,zalinescu2002convex}.

We will make use of \eqref{problem:uq_sdp_recovery} to recover an optimal probability measure for \eqref{problem:uq} under suitable sign assumptions. The relationship between these problems is depicted in Figure \ref{fig:recovery}. 

\usetikzlibrary{shapes.geometric, arrows}
\tikzstyle{process} = [rectangle, minimum width=3cm, minimum height=1cm, text centered, draw=black, fill=blue!10]
\tikzstyle{arrow} = [thick,->,>=stealth]
\begin{tikzpicture}[node distance=2cm]
<TikZ code>
\end{tikzpicture}

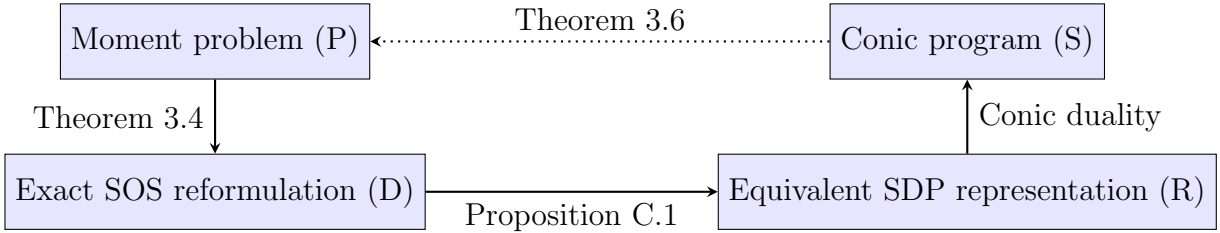
\begin{figure}[H]
    \centering
\begin{tikzpicture}[node distance=2cm]
\node (p) [process] {Moment problem \eqref{problem:uq}};

\node (d) [process, below of=p] {Exact SOS reformulation \eqref{problem:uq_dual}};

\node (sdp) [process, right of=d, xshift=8cm] {Equivalent SDP representation \eqref{problem:uq_sdp}};

\node (dd) [process, right of=p, xshift=8cm] {Conic program \eqref{problem:uq_sdp_recovery}};

\draw [arrow] (p) -- node[anchor=east] {Theorem \ref{thm:uqd_max_lin}} (d);

\draw[arrow] (d) -- node[anchor=north] {Proposition \ref{prop:sos-sdp} } (sdp);

\draw[arrow] (sdp) -- node[anchor=west] {Conic duality} (dd);

\draw[dotted,arrow] (dd) -- node[anchor=south] {Theorem \ref{thm:uq_recovery}} (p);
\end{tikzpicture}

\caption{Recovering an optimal probability measure for \eqref{problem:uq} from \eqref{problem:uq_sdp_recovery}. } \label{fig:recovery}
\end{figure}

We make use of the following generalized Jensen's inequality \cite[Theorem~2.6]{lasserre2009convexity}. 

\begin{proposition} \label{prop:jensens}
    Let $f \in \R[v]$ be an SOS-convex polynomial of an even degree $d$, $y = (y_\alpha)_{\alpha \in \mathcal{N}}$ satisfy $y_0 = 1$ and $\displaystyle \sum_{\alpha \in \mathcal{N}} y_\alpha B_\alpha \succeq 0$, and $L_y : \R[v] \to \R$ be the linear functional $\displaystyle L_y(f) = \sum_{\alpha \in \mathcal{N}} (f)_\alpha y_\alpha$. Then, $L_y (f) \geq f(L_y (\mathbf{v}))$, where $L_y (\mathbf{v}) = (L_y (\mathbf{v}_1), \ldots, L_y (\mathbf{v}_m))$, and $\mathbf{v}_i$ denotes the polynomial that maps a vector $v\in \R^m$ to its $i$-th coordinate $v_i$, $i=1,\ldots,m$.
\end{proposition}

\begin{theorem}
[\textbf{Recovering an optimal solution}] \label{thm:uq_recovery}
    Suppose that \eqref{problem:uq} admits a minimizer, Assumptions \ref{assump:interior}-\ref{assump:uq} are satisfied, and strong duality between \eqref{problem:uq_sdp} and \eqref{problem:uq_sdp_recovery} holds, i.e., $\max \eqref{problem:uq_sdp} = \min \eqref{problem:uq_sdp_recovery}$. Let $(\bar y^k, \bar \xi^k, \bar z^k) \in \R^{s(m,d)} \times \R^N \times \R$, $\bar y^k = (\bar y^k_\alpha)_{\alpha \in \mathcal{N}}$, be a minimizer for \eqref{problem:uq_sdp_recovery} with $\bar y^k_0 \geq 0$ for all $k=1,\ldots,r$. Denote $K := \{k \in \{1, \ldots,r \}: \overline{y}_0^k > 0\}$ and $\widehat{u}^k := \frac{1}{\bar y_0^k} (\bar y_{e_1}^k, \ldots, \bar y_{e_m}^k) \in \R^m$ for $k \in K$. Suppose that, for all $k \notin K$,  \begin{align*} 
        \overline{z}_k \ge 0 \mbox{ and } \sum_{\alpha \in \mathcal{N}} \overline{y}_{\alpha}^k (h_j)_{\alpha} \ge 0 \mbox{ for all } j=1,\ldots,J.
    \end{align*}
    
Then, $\displaystyle \widehat\mu := \sum_{k \in K} \bar y_0^k \1_{\widehat{u}^k}$, the linear combination of Dirac measures at points $\widehat{u}^k \in \R^m$, $k\in K$, is a minimizer for \eqref{problem:uq}.    
\end{theorem}

\begin{proof}
    The equalities $\min \eqref{problem:uq} = \max \eqref{problem:uq_dual}$, $\max \eqref{problem:uq_dual} = \max \eqref{problem:uq_sdp}$, and $\max \eqref{problem:uq_sdp} = \min \eqref{problem:uq_sdp_recovery}$ follow from Theorem \ref{thm:uqd_max_lin}, Proposition \ref{prop:sos-sdp}, and the strong duality assumption, respectively. Hence, $\min \eqref{problem:uq} = \min \eqref{problem:uq_sdp_recovery}$. Let $(\bar{y}^k, \bar \xi^k, \bar z^k) \in \mathbb{R}^{s(m,d)} \times \R^N \times \R$, $k=1,\ldots,r$, be a solution for problem \eqref{problem:uq_sdp_recovery} with $\bar y^k_0 \geq 0$ for all $k=1,\ldots,r$, so it satisfies Eqns. \eqref{eqn:recovery1}-\eqref{eqn:recovery4}. 

    For $k \in K$, define $\widehat{y}^k := \frac{1}{\bar y_0^k} \bar y^k \in \R^{s(m,d)}$, which satisfies $\widehat{y}_0^k = 1$ and by Eqn. \eqref{eqn:recovery4}, \begin{align*}
        \sum_{\alpha \in \mathcal{N}} \widehat{y}_\alpha^k B_\alpha = \frac{1}{\bar y_0^k} \sum_{\alpha\in \mathcal{N}} \bar y_\alpha^k B_\alpha \succeq 0.
    \end{align*}

    Let $\widehat u^k = (\widehat u_1^k, \ldots, \widehat u_m^k) = (\widehat y_{e_1}^k, \ldots, \widehat y_{e_m}^k)$, $k\in K$. Clearly, $\widehat{u}^k \in \Omega$ for all $k \in K$ because we can take $\widehat {\xi}^k := \frac{1}{\bar y_0^k} \bar \xi^k \in \R^N$ such that, by Eqn. \eqref{eqn:recovery3}, \begin{align*}
        F_0 + \sum_{i=1}^m \widehat{u}^k_i F_i + \sum_{t=1}^N \widehat {\xi}_t^k M_t = \frac{1}{\bar y_0^k} \Big (\bar y_0^k F_0 + \sum_{i=1}^m \bar y_{e_i}^k F_i + \sum_{t=1}^N \bar \xi_t^k M_t \Big ) \succeq 0,\ k \in K.
    \end{align*} 

    Define $\widehat \mu := \sum_{k \in K} \bar y_0^k \1_{\widehat u^k}$. We note that for any $f \in X'$, \begin{align*}
        \E_{\widehat\mu}^\Omega [f (\omega)] = \langle f, \widehat \mu \rangle = \langle f, \sum_{k\in K} \bar y_0^k \1_{\widehat u^k} \rangle = \sum_{k\in K} \bar y_0^k \langle f, \1_{\widehat u^k} \rangle = \sum_{k\in K} \bar y_0^k f (\widehat u^k). 
    \end{align*}

    By Eqn. \eqref{eqn:recovery4}, $\widehat \mu$ is a probability measure supported on $\Omega$ since $\E_{\widehat \mu}^\Omega [1] = \sum_{k\in K} \bar y_0^k = \sum_{k=1}^r \bar y_0^k = 1$.

    Now, $\widehat{u}^{k} = L_{\widehat{y}^{k}} (\mathbf{v})$, giving $h_j (\widehat u^k) = h_j (L_{\widehat y^k} (\mathbf{v}) )$ for $j=1,\ldots,J$, $k\in K$. Using Proposition \ref{prop:jensens} and the fact that $h_j$'s are SOS-convex polynomials, we have \begin{align*}
        h_j (\widehat u^k) \leq L_{\widehat y^k} (h_j) = \sum_{\alpha \in \mathcal{N}} \widehat y_\alpha^k (h_j)_\alpha = \frac{1}{\bar y_0^k} \sum_{\alpha\in \mathcal{N}} \bar y_\alpha^k (h_j)_\alpha, \ j=1,\ldots,J,\ k\in K. 
    \end{align*} 

    By the assumption that $\displaystyle \sum_{k \notin K} \sum_{\alpha \in \mathcal{N}} \bar y_\alpha^k (h_j)_\alpha \geq 0$ and Eqn. \eqref{eqn:recovery1}, it follows that \begin{align*}
        \E_{\widehat \mu}^\Omega [h_j (\omega)] = \sum_{k\in K} \bar y_0^k h_j (\widehat u^k) \leq \sum_{k\in K} \sum_{\alpha \in \mathcal{N}} \bar y_\alpha^k (h_j)_\alpha \leq \sum_{k=1}^r \sum_{\alpha \in \mathcal{N}} \bar y_\alpha^k (h_j)_\alpha \leq \gamma_j,\ j=1,\ldots,J,
    \end{align*}so $\widehat\mu$ is feasible for \eqref{problem:uq}. 

    Similarly, one can derive that $\displaystyle \bar y_0^k g_\ell ^k (\widehat u^k) \leq \sum_{\alpha \in \mathcal{N}} \bar y_\alpha^k (g_\ell^k)_\alpha$, $\ell=1,\ldots,L$, $k\in K$. From Eqn. \eqref{eqn:recovery2}, for each $k\in K$, $\displaystyle \sum_{\alpha \in \mathcal{N}} \bar y_\alpha^k (g_\ell^k)_\alpha \leq \bar z^k$. This shows that $\bar y_0^k g_\ell^k (\widehat u^k) \leq \bar z^k$ for all $\ell = 1,\ldots, L$, and so, $\displaystyle \max_{\ell=1,\ldots,L} \bar y_0^k g_\ell^k (\widehat u^k) \leq \bar z^k$. Let $\displaystyle g(v)= \min_{k=1,\ldots,r} \max_{\ell=1,\ldots,L} g_\ell^k(v)$. Recall that $g \in X'$ and, for each $k=1,\ldots,r$,  $g(v) \le \displaystyle \max_{\ell=1,\ldots,L} g_\ell^k(v)$ for any $v \in \Omega$. Then, since $\bar y_0^k>0$ for all $k \in K$, one has \begin{align*}
        \E_{\widehat \mu}^\Omega  [g (\omega) ] = \sum_{k \in K} \bar y_0^k \, g (\widehat u^k) \leq \sum_{k \in K} \bar y_0^k \, \left ( \max_{\ell=1,\ldots,L} g_\ell^k (\widehat u^k) \right ) = \sum_{k \in K}  \, \left ( \max_{\ell=1,\ldots,L} \bar y_0^k g_\ell^k (\widehat u^k) \right ) \leq \sum_{k\in K} \bar z^k. 
    \end{align*}
    
    Using the assumption that $\displaystyle \sum_{k\notin K} \bar z^k \geq 0$, one has $\displaystyle \E_{\widehat\mu}^\Omega [g(\omega)] \leq \sum_{k\in K} \bar z^k \leq \sum_{k=1}^r \bar z^k  = \min \eqref{problem:uq_sdp_recovery}$, giving $\min\eqref{problem:uq} \leq \min \eqref{problem:uq_sdp_recovery}$.   
    
    Conversely, $\min \eqref{problem:uq} = \min \eqref{problem:uq_sdp_recovery}$ ensures that $\displaystyle \min\eqref{problem:uq} = \E_{\widehat \mu}^\Omega \Big [\min_{k=1,\ldots,r} \max_{\ell=1,\ldots,L} g_\ell^k (\omega) \Big ]$. Therefore, $\widehat \mu$ is an optimal solution for the original moment problem \eqref{problem:uq}. 
\end{proof}

The optimal probability measure $\widehat \mu := \sum_{k\in K} \bar y_0^k \1_{\widehat{u}^k}$ refers to the discrete distribution with probabilities $\mathbb{P} (\omega = \widehat u^k) = \bar y_0^k$, $k\in K$, i.e., the distribution with masses of $\bar y_0^k$ at points $\widehat u^k \in \Omega$, $k\in K$. 

We present a recovery result under a stronger condition that is easy to verify. 

\begin{corollary}[\textbf{Recovery under a simple stronger condition}] \label{coro:recovery2}
     Suppose that \eqref{problem:uq} admits a minimizer, Assumptions \ref{assump:interior}-\ref{assump:uq} are satisfied, and $\max \eqref{problem:uq_sdp} = \min \eqref{problem:uq_sdp_recovery}$. Let $( \bar y^k, \bar \xi^k, \bar z^k) \in \R^{s(m,d)} \times \R^N \times \R$, $\bar y^k = (\bar y^k_\alpha)_{\alpha \in \mathcal{N}}$, be a minimizer for \eqref{problem:uq_sdp_recovery} with $\bar y^k_0 > 0$ for all $k=1,\ldots,r$. Denote $\widehat{u}^{k} := \frac{1}{\bar y_0^k}(\bar{y}_{e_1}^k,\ldots,\bar{y}_{e_m}^k) \in \R^m$, $k=1,\ldots,r$. Then, $\widehat\mu := \sum_{k=1}^r \bar y_0^k \1_{\widehat{u}^k}$ is a minimizer for \eqref{problem:uq}.  
\end{corollary}

\begin{proof}
    This follows from Theorem \ref{thm:uq_recovery} by setting $K = \{1,\ldots,r\}$. 
\end{proof}

As an application to polynomial optimization, our recovery result reduces to the known Lasserre's result in  \cite{lasserre2009convexity}. More explicitly, consider \begin{align} \label{problem:poly}
    \min_{v\in \Omega} g (v),
\end{align}where $\Omega\subset \R^m$ is given as in \eqref{eqn:lmi} and $g$ is an SOS-convex polynomial. Notice that, in \cite[Theorem~3.3]{lasserre2009convexity}, the feasible set $\Omega$ is a compact semi-algebraic set formed by SOS-convex polynomial inequalities, and thus can be represented as a projected spectrahedron \cite{helton2008structured}.

We associate \eqref{problem:poly} with the following conic programs  \begin{align} \label{problem:poly_sdp}
    \max_{\lambda \in \R, Z \in \mathbb S_+^\nu} \Big \{ -\lambda : g (v) + \lambda - \tr(ZF_0) - \sum_{i=1}^m v_i \tr (Z F_i) \in \Sigma_d^2 (v),\ \tr(ZM_t) = 0,\ t=1,\ldots,N \Big \},
\end{align}and its dual problem \begin{align} \label{problem:poly_recovery}
    \min_{y \in \R^{s(m,d)}, \xi \in \R^N} \Big \{ \sum_{\alpha \in \mathcal{N}} y_\alpha (g)_\alpha : F_0 + \sum_{i=1}^m y_{e_i} F_i + \sum_{t=1}^N \xi_t M_t \succeq 0,\ y_0 = 1,\ \sum_{\alpha \in \mathcal{N}} y_\alpha B_\alpha \succeq 0 \Big \}.
\end{align}

The following Corollary shows how an optimal solution of \eqref{problem:poly} can be obtained by solving the SDP program \eqref{problem:poly_recovery}. 

\begin{corollary} For the problem \eqref{problem:poly}, 
    suppose that $\Omega$ is compact, there exist $\bar v \in \R^m$ and $\bar \xi \in \R^N$ such that $F_0 + \sum_{i=1}^m \bar v_i F_i + \sum_{t=1}^N \bar \xi_t M_t \succ 0$, and $\max \eqref{problem:poly_sdp} = \min \eqref{problem:poly_recovery}$. Let $(\bar y, \bar \xi) \in \R^{s(m,d)} \times \R^N$ be a minimizer for \eqref{problem:poly_recovery}. Then, $(\bar y_{e_1}, \ldots, \bar y_{e_m}) \in \R^m$ is a minimizer for \eqref{problem:poly}. 
\end{corollary}

\begin{proof}
    Consider the moment problems \begin{align} \label{eqn:poly_dirac0}
        \min_{\mu \in \mathcal{P}_\Omega} \ \E_\mu^\Omega [g (\omega)]
    \end{align}and \begin{align} \label{eqn:poly_dirac}
        \min_{\mu \in D} \ \E_\mu^\Omega [g (\omega)],
    \end{align}where $D = \{ \1_v \in X : v \in \Omega \} \subset \mathcal{P}_\Omega$ is the set of all Dirac measures supported on $\Omega$, giving $\min \eqref{eqn:poly_dirac0} \leq \min \eqref{eqn:poly_dirac}$. Note that $\min \eqref{eqn:poly_dirac} = \min \eqref{problem:poly}$ because $\E_{\1_v}^\Omega [g(\omega)] = g(v)$ for all $v\in \Omega$. 
    
    The interior point condition (Assumption \ref{assump:interior}) for the moment problem \eqref{eqn:poly_dirac0} requires that $1 \in \interior \{ \E_\mu^\Omega [1] : \mu \in C\}$, where $C = \{\mu \in X : \mu \succeq_\B 0\}$. Since $C = \text{cone}(\mathcal{P}_\Omega)$ and $\E_{\mathbb P}^\Omega [1] = 1$ for any $\mathbb P \in \mathcal{P}_\Omega$, a direct verification shows that 
        $\{\E_\mu^\Omega [1] : \mu \in C\} 
        = \R_+$. So, Assumption \ref{assump:interior} is satisfied automatically. 
    
    By Corollary \ref{coro:recovery2} with $\bar y_0 = 1$, we have that $\widehat \mu := \1_{(\bar y_{e_1}, \ldots, \bar y_{e_m})}$ is a minimizer for \eqref{eqn:poly_dirac0} and $(\bar y_{e_1}, \ldots, \bar y_{e_m}) \in \Omega$. As $\E_{\widehat{\mu}}^\Omega [g (\omega)]=\min \eqref{eqn:poly_dirac0} \leq \min \eqref{eqn:poly_dirac}$ and $\widehat \mu \in D$, so $\widehat \mu$ is also a minimizer for \eqref{eqn:poly_dirac}. Noting that $\E_{\widehat \mu }^\Omega [g (\omega)] = g (\bar y_{e_1}, \ldots, \bar y_{e_m})$, one conclude that $(\bar y_{e_1}, \ldots, \bar y_{e_m})$ is a minimizer for the polynomial optimization problem \eqref{problem:poly}. 
\end{proof}

\section{Applications to Newsvendor and Revenue Maximization} \label{sec:applications} 

We devote this section to presenting SDP reformulations and numerical results for the Newsvendor and revenue maximization problems. 

\textbf{Generalized Newsvendor Problems}. A company orders $n$ goods, where $n \geq 1$. For each of the $i$-th goods, the unit cost of upfront ordering is $\$c_i$ with $0<c _ i<1$. For a fixed order quantity $0 \leq x_i \leq R_i$, where $R_i$ refers to the capacity, the company wishes to estimate a worst-case (upper) bound for the cost of ordering to meet the random demands. If the demand $\omega_i$ exceeds the upfront order $x_i$, the company incurs a unit back-ordering cost of \$1, or a total back-ordering cost of $\$ (\omega_i - x_i)$, otherwise, the cost of back-ordering is $\$ 0$. This results in a total ordering cost of $c_i x_i + \max \{\omega_i - x_i, 0 \}$. 

Assume that the demand for each goods is independent and the support is $\prod_{i=1}^n \Omega_i$, where $\Omega_i = [\underline \omega_i, \overline \omega_i]$ with $0 \leq \underline \omega_i < \overline \omega_i$, $i=1,\ldots,n$. Then, the multi-product Newsvendor model reduces to single-product Newsvendor models, and the worst-case ordering cost can be found by solving $\displaystyle \max_{\mu \in \probset_i} \{ c_i x_i + \E_\mu^{\Omega_i} [ \max\{ \omega_i - x_i ,0 \} ] \}$ for all $i=1,\ldots,n$, where $\probset_i$ consists of various moment constraints. 

We focus on the single-product Newsvendor model and suppress the index $i \in \{1,\ldots,n\}$. For a fixed order quantity $0\leq x\leq R$, the worst-case ordering cost can be found from the following moment problem, \begin{align} \label{problem:Newsvendor_uq}
   \max_{\mu \in \probset}\ cx + \E_\mu^\Omega [\max \{\omega - x, 0\}],
\end{align}where $\Omega = [\underline  \omega, \overline\omega] = \{v\in \R : F_0 + v F_1 \succeq 0\}, \underline \omega < \overline \omega$, with $F_0 = \begin{pmatrix}
    -\underline \omega & 0 \\ 0 & \overline \omega
\end{pmatrix}$ and $F_1 = \begin{pmatrix}
    1 & 0 \\ 0 & -1
\end{pmatrix}$.

Another interpretation of \eqref{problem:Newsvendor_uq} at $c=0$ comes from pricing a European call option \cite{bertsimas2002relation,lasserre2008semidefinite}. The random variable $\omega$ represents the price of the underlying stock, and the fixed value $x$ is regarded as the strike price. This gives the expected price of $\E_\mu^\Omega [\max \{\omega - x, 0\}]$, and thus a sharp upper bound $\displaystyle \max_{\mu \in \probset} \E_\mu^\Omega [\max\{\omega- x,0\}]$ of the price. 

Numerically tractable reformulations for \eqref{problem:Newsvendor_uq} are available for a limited choice of $\probset$. For example,  \cite{lo1987semi} provides a closed-form formula when $\probset$ specifies known values for the mean and variance; recently, \cite{guo2022unified} presents a closed-form formula when $\probset$ specifies known values for mean plus one other moment of any order; \cite{bertsimas2002relation} reformulates \eqref{problem:Newsvendor_uq} as an SDP when $\probset$ contains the first $J$ moments; and \cite{wiesemann2014distributionally} reformulates \eqref{problem:Newsvendor_uq} as a conic program when the random demand variable follows a distribution defined by robust statistics. For other related models with tractable reformulations, see  \cite{delage2010distributionally}.  

We study the problem \eqref{problem:Newsvendor_uq} when $\probset_1 = \{\mu \in \mathcal{P}_\Omega : \E_\mu^\Omega [\omega] \leq \gamma_1, \E_\mu^\Omega [\omega^2] \leq \gamma_2\}$ and $\probset_2 = \probset_1 \cap \{ \mu \in \mathcal{P}_\Omega : \E_\mu^\Omega [\omega^4]\leq \gamma_3\}$, where $\mathcal{P}_\Omega$ is the set of probability measures on $\Omega$. The set $\probset_1$ contains upper bounds for the mean and the second-order moment which is related to the variance, while $\probset_2$ contains an additional bound for the fourth-order moment which is related to the kurtosis of a distribution. The mean, variance, and kurtosis of the random demand distribution are usually accessible in real-world applications. 

The standard minimization form for the Newsvendor problem is \begin{align} \label{problem:newsvendor_std}
    \min \eqref{problem:newsvendor_std} := \min_{\mu \in \probset}\ \E_\mu^\Omega \Big [\min_{k=1,2} g_1^k (\omega) \Big ], \ g_1^1 (v) = x- v,\ g_1^2 (v) = 0,
\end{align}and the worst-case expected ordering cost is equal to $[cx - \min \eqref{problem:newsvendor_std}]$.

The associated SOS optimization problem for  \eqref{problem:newsvendor_std} when $\probset = \probset_1$ is  \begin{align} \label{problem:newsvendor_sos1}
    \max \eqref{problem:newsvendor_sos1} := \max_{ \substack{ \lambda \in \R_+^2\times \R \\ Z, Z' \in \mathbb S_+^2} }\ & - \lambda_1 \gamma_1 - \lambda_2 \gamma_2 - \lambda_3 \\ \st \quad & [x + \lambda_3 - \tr (ZF_0) ]+ [-1 + \lambda_1 - \tr (ZF_1)] v + \lambda_2 v^2 \in \Sigma_2^2 (v), \nonumber \\ & [\lambda_3 - \tr(Z' F_0)] + [\lambda_1 - \tr(Z' F_1)] v + \lambda_2 v^2 \in \Sigma_2^2 (v). \nonumber
\end{align}

Its equivalent SDP representation is \begin{align} \label{problem:Newsvendor_sdp1} 
    \max \eqref{problem:Newsvendor_sdp1} := \max_{\substack {\lambda \in \R_+^2 \times \R\\ Z, Z' \in \mathbb S_+^2} }\ & - \lambda_1 \gamma_1  - \lambda_2 \gamma_2 - \lambda_3 \\ \st \quad   & \begin{pmatrix}        x + \lambda_3 - \tr(ZF_0) & \frac{1}{2} [-1 + \lambda_1  -\tr(Z F_1) ] \\ \frac{1}{2} [-1 + \lambda_1 - \tr (Z F_1) ] & \lambda_2    \end{pmatrix} \succeq 0, \nonumber \\ & \begin{pmatrix}        \lambda_3 - \tr(Z' F_0) & \frac{1}{2} [\lambda_1 - \tr(Z' F_1)] \\ \frac{1}{2} [\lambda_1 - \tr(Z'F_1)] & \lambda_2      \end{pmatrix} \succeq 0. \nonumber
\end{align}

\begin{proposition} \label{prop:newsvendor_uq1}
    Assume that \eqref{problem:newsvendor_std} admits a minimizer for $\probset = \probset_1$ and $(\gamma_1 , \gamma_2 , 1 ) \in$ 
    
    $\interior \{ \{ (\E_\mu^\Omega [\omega], \E_\mu^\Omega [\omega^2], \E_\mu^\Omega [1] ) : \mu \in X, \mu \succeq_\B 0 \} + \R_+^2 \times \{0\} \}$. Then, $\min \eqref{problem:newsvendor_std} =\max \eqref{problem:Newsvendor_sdp1}$. 
\end{proposition}

\begin{proof} 
    By Theorem \ref{thm:uqd_max_lin}, $\min \eqref{problem:newsvendor_std} = \max \eqref{problem:newsvendor_sos1}$. Next, following the approach in Proposition \ref{prop:sos-sdp}, the constraint $[x + \lambda_3 - \tr(ZF_0)] + v[-1 + \lambda_1 - \tr(ZF_1)] + \lambda_2 v^2 \in \Sigma_2^2 (v)$ is equivalent to the existence of $Q \in \mathbb S_+^2$ such that $\tr(QB_0) = x + \lambda_3 - \tr(ZF_0)$, $\tr(QB_1) = -1 + \lambda_1 - \tr(ZF_1)$, and $\tr(QB_2) = \lambda_2$, where $B_0 = \begin{pmatrix}
        1 & 0 \\ 0 & 0
    \end{pmatrix}$, $B_1 = \begin{pmatrix}
        0 & 1 \\ 1 & 0
    \end{pmatrix}$, and $B_2 = \begin{pmatrix}
        0 & 0 \\ 0 & 1
    \end{pmatrix}$. This gives \begin{align*}
        Q = \begin{pmatrix}        x + \lambda_3 - \tr(ZF_0) & \frac{1}{2} [-1 + \lambda_1  -\tr(Z F_1) ] \\ \frac{1}{2} [-1 + \lambda_1 - \tr (Z F_1) ] & \lambda_2    \end{pmatrix} \in \mathbb S_+^2. 
    \end{align*}
    
    Similar arguments apply for the other constraint, and hence $\max \eqref{problem:newsvendor_sos1} = \max \eqref{problem:Newsvendor_sdp1}$. 
\end{proof}

Now, we illustrate how one can recover an optimal measure for \eqref{problem:newsvendor_std} via the following SDP \begin{align} \label{problem:newsvendor_recovery1}
    \min \eqref{problem:newsvendor_recovery1} := \min_{\substack{y^1, y^2 \in \R^3 \\ z^1, z^2 \in \R}}\ &  z^1 + z^2 \\ \st \quad & y_1^1 + y_1^2 \leq \gamma_1,\ y_2^1 + y_2^2 \leq \gamma_2,\ y_0^1 + y_0^2 = 1, \nonumber \\ & x y_0^1 - y_1^1 \leq z^1,\ 0 \leq z^2 , \nonumber \\ & y_0^k F_0 + y_1^k F_1 \succeq 0,\ \begin{pmatrix}
        y_0^k & y_1^k \\ y_1^k & y_2^k
    \end{pmatrix} \succeq 0, \ k=1,2. \nonumber 
\end{align}

\begin{proposition}
    Assume the same conditions as in Proposition \ref{prop:newsvendor_uq1} and $\max \eqref{problem:Newsvendor_sdp1} = \min \eqref{problem:newsvendor_recovery1}$. Let $(\bar y^k, \bar z^k) \in \R^3\times \R$ with $\bar y^k=(\bar y^k_0, \bar y_1^k, \bar y_2^k)$ be a minimizer for \eqref{problem:newsvendor_recovery1} with $\bar y^k_0 > 0$, $k=1,2$. Denote $\widehat{u}^k := \frac{1}{\bar y_0^k} \bar y_1^k \in \R$, $k=1,2$. Then, $\bar y_0^1 \1_{\widehat{u}^1} + \bar y_0^2 \1_{\widehat{u}^2}$ is a minimizer for \eqref{problem:newsvendor_std} for $\probset = \probset_1$. 
\end{proposition}

\begin{proof}
    This follows from Corollary \ref{coro:recovery2}. 
\end{proof}

Note that the probability measure optimal for the standard minimization form \eqref{problem:newsvendor_std} will also be optimal for the original Newsvendor maximization problem \eqref{problem:Newsvendor_uq}. 

In a similar vein, we examine the Newsvendor model with $\probset = \probset_2$. Its associated SOS optimization problem is given by \begin{align} \label{problem:newsvendor_sos2}
    \max \eqref{problem:newsvendor_sos2} := \max_{\substack{ \lambda\in \R^3_+\times \R \\ Z, Z' \in \mathbb S_+^2} } \ & - \lambda_1 \gamma_1 - \lambda_2 \gamma_2 - \lambda_3 \gamma_3 - \lambda_4 \\ \st \quad & [ x + \lambda_4 - \tr(Z F_0) ] + [-1 + \lambda_1 - \tr(Z F_1) ] v + \lambda_2 v^2 + \lambda_3 v^4 \in \Sigma_4^2 (v), \nonumber \\ & [\lambda_4 - \tr(Z' F_0) ] + [\lambda_1 - \tr(Z' F_1) ] v + \lambda_2 v^2 + \lambda_3 v^4 \in \Sigma_4^2(v). \nonumber
\end{align}

Its equivalent SDP representation is \begin{align} \label{problem:newsvendor_sdp2}
    \max \eqref{problem:newsvendor_sdp2}  := \max_{\substack {\lambda \in \R^3_+ \times \R \\  Z, Z' \in \mathbb S_+^2, Q, Q' \in \mathbb S_+^3 } }\ & - \lambda_1 \gamma_1 - \lambda_2 \gamma_2 - \lambda_3 \gamma_3 - \lambda_4 \\ \st \quad\quad\quad & \tr(QB_0) = x + \lambda_4 - \tr(Z F_0),\ \tr(Q' B_0) = \lambda_4 - \tr(Z' F_0), \nonumber \\ & \tr (Q B_1) = -1+ \lambda_1 - \tr(ZF_1),\ \tr( Q' B_1) = \lambda_1 - \tr(Z' F_1), \nonumber \\ & \tr(Q B_2) = \tr(Q' B_2) = \lambda_2, \nonumber \\ & \tr(Q B_3) = \tr (Q' B_3) = 0, \nonumber \\ & \tr( QB_4) = \tr(Q' B_4) = \lambda_3, \nonumber
\end{align}where \begin{align} \label{eqn:B03}
    B_0 = \begin{pmatrix}
        1 & 0 & 0 \\ 0 & 0 & 0 \\ 0 & 0 & 0
    \end{pmatrix}, B_1 = \begin{pmatrix}
        0 & 1 & 0 \\ 1 & 0 & 0 \\ 0 & 0 & 0
    \end{pmatrix}, B_2 = \begin{pmatrix}
        0& 0 & 1 \\ 0 & 1 & 0 \\ 1 & 0 & 0 
    \end{pmatrix}, B_3 = \begin{pmatrix}
        0 & 0 & 0 \\ 0 & 0 & 1 \\ 0 & 1 & 0
    \end{pmatrix}, B_4 = \begin{pmatrix}
        0 & 0 & 0 \\ 0 & 0 & 0 \\ 0 & 0 & 1
    \end{pmatrix}.
\end{align}

\begin{proposition} \label{prop:newsvendor_uq2}
    Assume that \eqref{problem:newsvendor_std} admits a minimizer for $\probset = \probset_2$ and $(\gamma_1, \gamma_2, \gamma_3,1 ) \in \interior \{ \{ (\E_\mu^\Omega [\omega], \E_\mu^\Omega [\omega^2], \E_\mu^\Omega [\omega^4], \E_\mu^\Omega [1] ) : \mu \in X, \mu \succeq_\B 0 \} + \R_+^3 \times \{0\} \}$. Then, $\min \eqref{problem:newsvendor_std} = \max \eqref{problem:newsvendor_sdp2}$. 
\end{proposition}

\begin{proof}
    The proof is similar to Proposition \ref{prop:newsvendor_uq1}. 
\end{proof}

We recover an optimal solution for \eqref{problem:newsvendor_std} using the following SDP, \begin{align} \label{problem:newsvendor_recovery2}
    \min \eqref{problem:newsvendor_recovery2} := \min_{\substack{y^1, y^2 \in \R^5 \\ z^1, z^2 \in \R}}\ & z^1 + z^2 \\ \st\quad \ & y^1_1 + y^2_1 \leq \gamma_1, \ y_2^1 + y_2^2 \leq \gamma_2,\ y_4^1+y_4^2 \leq \gamma_3 ,\ y_0^1 + y_0^2 = 1, \nonumber \\ & x y_0^1 - y_1^1 \leq z^1, \ 0\leq z^2 , \nonumber \\ & y_0^k F_0 + y_1^k F_1 \succeq 0,\ \begin{pmatrix}
        y_0^k & y_1^k & y_2^k \\ y_1^k & y_2^k & y_3^k \\ y_2^k & y_3^k & y_4^k
    \end{pmatrix}  \succeq 0,\ k=1,2. \nonumber
\end{align}

\begin{proposition}
    Assume the same conditions as in Proposition \ref{prop:newsvendor_uq2} and $\max \eqref{problem:newsvendor_sdp2} = \min \eqref{problem:newsvendor_recovery2}$. Let $(\bar y^k, \bar z^k) \in \R^5 \times \R$ with $\bar y^k = (\bar y_0^k, \bar y_1^k, \bar y_2^k, \bar y_3^k, \bar y_4^k)$ be a minimizer for \eqref{problem:newsvendor_recovery2} with $\bar y_0^k > 0$, $k=1,2$. Denote $\widehat u^k := \frac{1}{\bar y_0^k} \bar y_1^k \in \R$, $k=1,2$. Then, $\bar y_0^1 \1_{\widehat u^1} + \bar y_0^2 \1_{\widehat{u}^2}$ is a minimizer for \eqref{problem:newsvendor_std} for $\probset = \probset_2$. 
\end{proposition}

\begin{proof}
    This follows from Corollary \ref{coro:recovery2}. 
\end{proof}

\begin{figure}[H]
    \centering
    \includegraphics[width=0.4\linewidth]{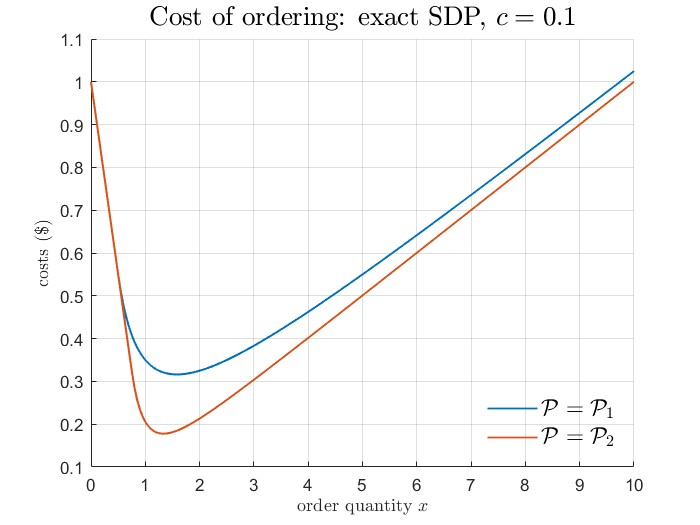}
    \caption{Costs of ordering (\$) for $\probset = \probset_1$ (blue) and $\probset = \probset_2$ (orange). }
    \label{fig:newsvendor}
\end{figure}

We illustrate how the worst-case costs of ordering $[cx-\max\eqref{problem:Newsvendor_sdp1}]$ and $[cx - \max\eqref{problem:newsvendor_sdp2}]$ change by varying the order quantity $x \in [0, R] = [0,10]$ for $\probset_1$ and $\probset_2$. The SDPs are modelled by CVX MATLAB \cite{cvx} and solved by Mosek \cite{mosek}. Set $\Omega = [\underline \omega, \overline \omega] = [0,100]$ for the support, $c = 0.1$ for the unit cost of upfront ordering, and $\gamma_1 = \gamma_2 = \gamma_3 = 1$ for the moment bounds. Figure \ref{fig:newsvendor} illustrates the cost of
ordering (\$) for $\probset = \probset_1$ (blue) and $\probset = \probset_2$ (orange). 

The curves in Figure \ref{fig:newsvendor} exhibit ``U'' shapes. Before the turning points ($x\in [ 0, 1.5811 ]$ blue, $x\in [0, 1.3337]$ orange), the more upfront orders at cost $c = \$0.1$ the company makes, the more it saves from back-ordering at cost \$1. For larger upfront ordering quantities ($x\in [1.5811, 10]$ blue, $x\in [1.3337, 10]$ orange) that can largely cover the demand, back-ordering is not needed. The slope of the curves in these regions is approximately $c = 0.1$.

Notice that $\probset_2$ is a subset of $\probset_1$, Figure \ref{fig:newsvendor} confirms that the costs of ordering for $\probset_1$ (blue) are at least the costs for $\probset_2$ (orange), and this leads to a higher minimum cost for $\probset_1$. This could be justified that, by knowing more information about the demand, i.e., $\probset = \probset_2$, the future demand can be better predicted, and thus the ordering strategy could be better determined. 

The minimum costs of ordering are achieved at the turning points of the curves. Specifically, the minimum costs are \$0.3162 and \$0.1778 if the company orders $x = 1.5811$ and $x = 1.3337$ units of goods under the ambiguity constraints $\probset = \probset_1$ (blue) and $\probset = \probset_2$ (orange), respectively. 

For $\probset_1$ with order quantity $x = 1.5811$, an optimal probability measure is $0.1 \1_{3.1623} + 0.9 \1_{0.0003}$. It refers to the discrete distribution with $\mathbb{P} (\omega = 3.1623) = 0.1$ and $\mathbb{P} (\omega = 0.0003) = 0.9$. For $\probset_2$ with $x = 1.3337$, an optimal probability measure is $0.1 \1_{1.7782} + 0.9 \1_{0.0200}$. 

\textbf{Revenue Maximization. }A merchant supplies a random quantity $\omega \in [0,R] \subseteq \R$ of goods and sells them to $n$ customers. Each customer offers a different price based on the supply quantity. The goods can be sold to one customer exclusively at a time, and the merchant wishes to maximize revenue by selling the goods to the customer who offers the highest price. 

This problem of revenue maximization can be formulated as \begin{align} \label{problem:revenue_uq}
    \max_{\mu \in \mathcal{P}_\Omega} \Big \{ \E_\mu^\Omega \Big [\max_{k=1,\ldots,n} h_k (\omega)  \Big ] : \E_\mu^\Omega [\omega] \leq \gamma_1, \E_\mu^\Omega [\omega^2] \leq \gamma_2 \Big \},
\end{align}where $h_k$ is the offer price from the $k$-th customer, $k=1,\ldots,n$, and $\Omega = [0,R]$, $R > 0$. The support can be expressed equivalently as a spectrahedron $\Omega = \{ v\in \R : F_0 + v F_1 \succeq 0\}$ with $F_0 = \begin{pmatrix}
    0 & 0 \\ 0 & R
\end{pmatrix}$ and $F_1 = \begin{pmatrix}
    1 & 0 \\ 0 & -1
\end{pmatrix}$.

The standard minimization form is \begin{align} \label{problem:revenue_std}
    \min \eqref{problem:revenue_std} := \min_{\mu \in \mathcal{P}_\Omega} \Big \{ \E_\mu^\Omega [ g (\omega) ] : \E_\mu^\Omega [\omega] \leq \gamma_1, \E_\mu^\Omega [\omega^2] \leq \gamma_2 \Big \},
\end{align}where $g$ is given by $\displaystyle g(v)= -\max_{k=1,\ldots,n} h_k (v)$. The maximum expected revenue is equal to $[-\min \eqref{problem:revenue_std}]$. Denoting $f_k(v)=-h_k(v)$, one has $g(v)=\displaystyle \min_{k=1,\ldots,n} f_k(v)$. Here, for each $k=1,\ldots,n$, we describe the offer price by \begin{align} \label{eqn:offer_price}    
    f_k(v) = & \begin{cases} \alpha_k (v - b_k)^2 + \beta_k (v-b_k)^4 + c_k,\ & \text{if}\ 0\leq v \leq b_k, \\ c_k,\ & \text{if}\ v > b_k,    \end{cases}
\end{align}where $\alpha_k, \beta_k, b_k \geq 0$ and $c_k \leq -\alpha_k b_k^2 - \beta_k b_k^4$. Note that $g$ is a piecewise SOS-convex function with the following representation: \begin{align} \label{problem:revenue_uq1}
    g(v)= \min_{k=1,\ldots,2n} \max_{\ell=1,2} g_\ell^k (v),
\end{align}where $g_1^k (v) = g_2^k (v) = \alpha_k (v- b_k)^2 + \beta_k (v-b_k)^4 + c_k$, $g_1^{n+k} (v) = -(\alpha_k b_k + \beta_k b_k^3) v + (\alpha_k b_k^2 + \beta_k b_k^4 + c_k)$, and $g_2^{n+k} (v) = c_k$, for $k=1,\ldots,n$. For details of this representation, see  Appendix \ref{sec:piecewise}. 

In fact, $(-f_k)$ where $f_k$ is given as in Eqn. \eqref{eqn:offer_price}, or $h_k$ in problem \eqref{problem:revenue_uq}, corresponds to a combined quadratic-quartic utility function that is concave, non-decreasing, and smooth \cite{gerber1998utility}. The coefficients $(-\alpha_k)$ and $(-\beta_k)$ capture the rate at which the customer increases the offer price based on the supply; $b_k$ is the maximum quantity of goods the $k$-th customer is willing to purchase, beyond which the customer will no longer be willing to increase the price; and $c_k$ which excludes negative offer prices is the maximum price. 

The function in Eqn. \eqref{eqn:offer_price} reduces to the quadratic model in \cite{han2015convex} by setting $\beta_k = 0$, $k=1,\ldots,n$, where an approximating upper bound for the maximum revenue is calculated through a convex program. Our model covers more diverse purchasing behaviours described by the quartic feature.

\begin{figure}[H]
    \centering
    \includegraphics[width=0.4\linewidth]{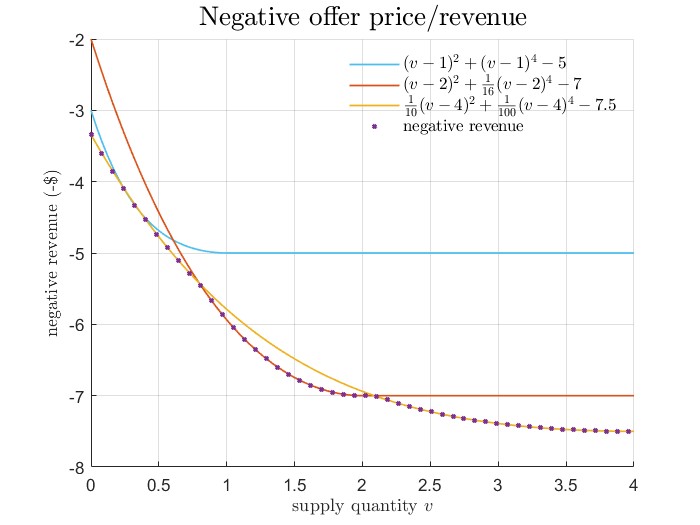}
    \includegraphics[width=0.4\linewidth]{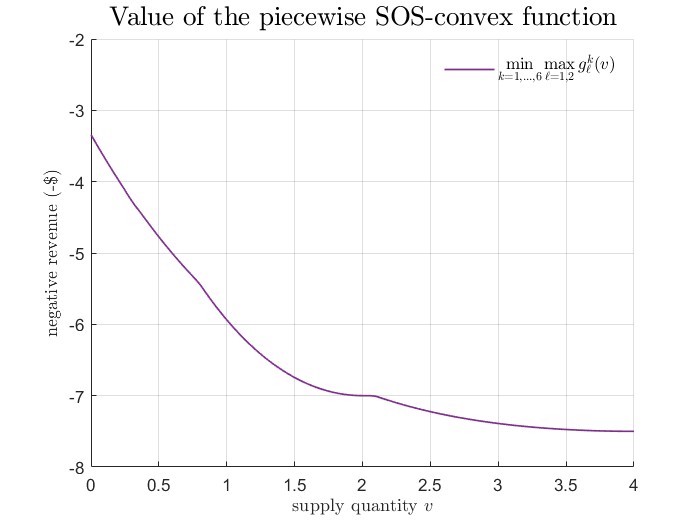}
    \caption{Left: negative offer prices (solid lines) and negative revenue (dotted lines). Right: piecewise SOS-convex function in Eqn. \eqref{problem:revenue_uq1}. }
    \label{fig:revenue}
\end{figure}

Suppose that there are three customers with parameters $\alpha_1 = 1$, $\alpha_2 = 1$, $\alpha_3 = \frac{1}{10}$, $\beta_1 = 1$, $\beta_2 = \frac{1}{16}$, $\beta_3 = \frac{1}{100}$, $b_1 = 1$, $b_2 = 2$, $b_3 = 4$, and $c_1 = -5$, $c_2 = -7$, $c_3 = -7.5$. The supply quantity ranges between $[0,R] = [0,4]$. As shown in Figure \ref{fig:revenue}, the piecewise SOS-convex function $g$ is neither convex nor concave and is not smooth.

Associate \eqref{problem:revenue_std} with the following SOS optimization problem with $n=3$ customers, \begin{align*} 
    \nonumber \\ \max_{\substack{\lambda \in \R_+^2 \times \R, \delta_k \in [0,1] \\ Z_k, Z_k' \in \mathbb S_+^2, k=1,\ldots,n}}\ & -\lambda_1 \gamma_1 - \lambda_2 \gamma_2 - \lambda_3 \\ \st \quad\quad \quad & [\alpha_k b_k^2 + \beta_k b_k^4 + c_k + \lambda_3 - \tr(Z_k F_0)] + [- 2 \alpha_k b_k - 4 \beta_k b_k^3 + \lambda_1 - \tr(Z_k F_1)] v +  \nonumber \\ &\quad \quad\quad [\alpha_k + 6\beta_k b_k^2 + \lambda_2] v^2 - 4 \beta_k b_k v^3 +  \beta_k v^4 \in \Sigma_4^2 (v),\ k=1,\ldots,n, \nonumber \\ & [\delta_k \alpha_k b_k^2 + \delta_k \beta_k b_k^4 + c_k + \lambda_3 - \tr(Z_k' F_0)] + [- \delta_k \alpha_k b_k - \delta_k \beta_k b_k^3 + \lambda_1 - \tr(Z_k' F_1) ] v + \nonumber \\ & \quad \quad \quad \lambda_2 v^2 \in \Sigma_4^2 (v),\ k=1,\ldots,n, \nonumber
\end{align*}and its equivalent SDP representation \begin{align} \label{problem:revenue_sdp1}
    \max \eqref{problem:revenue_sdp1} := 
    \max_{\substack{\lambda \in \R_+^2 \times \R, \delta_k \in [0,1] \\ Z_k, Z_k' \in \mathbb S_+^2, Q_k, Q_k' \in \mathbb S_+^3, k=1,\ldots,n}}\ & -\lambda_1 \gamma_1 - \lambda_2 \gamma_2  - \lambda_3 \\ \st \quad\quad\quad\quad & \text{for each}\ k=1,\ldots,n, \nonumber \\ & \tr(Q_k B_0) = \alpha_k b_k^2 + \beta_k b_k^4 + c_k + \lambda_3 - \tr(Z_k F_0), \nonumber \\ & \tr(Q_k B_1) = - 2 \alpha_k b_k - 4 \beta_k b_k^3 + \lambda_1 - \tr(Z_k F_1), \nonumber \\ & \tr(Q_k B_2) = \alpha_k + 6 \beta_k b_k^2 + \lambda_2, \nonumber \\ & \tr(Q_k B_3) = - 4 \beta_k b_k, \ \tr(Q_k B_4) = \beta_k, \nonumber \\ & \tr(Q_k' B_0) = \delta_k \alpha_k b_k^2 + \delta_k \beta_k b_k^4 + c_k + \lambda_3 - \tr(Z_k' F_0), \nonumber \\ & \tr(Q_k' B_1) = -\delta_k \alpha_k b_k - \delta_k \beta_k b_k^3 + \lambda_1 - \tr(Z_k' F_1), \nonumber \\ & \tr(Q_k' B_2) = \lambda_2,\ \tr(Q_k' B_3) = \tr(Q_k' B_4) = 0, \nonumber 
\end{align}where $B_0, \ldots, B_4$ are given in Eqn. \eqref{eqn:B03}.

\begin{proposition} \label{prop:revenue1}
    Assume that \eqref{problem:revenue_std} admits a minimizer and $(\gamma_1 , \gamma_2, 1) \in \interior \{ \{ (\E_\mu^\Omega [\omega] , \E_\mu^\Omega [\omega^2] , \E_\mu^\Omega [1]) : \mu \in X, \mu \succeq_\B 0 \} + \R_+^2 \times \{0\} \}$. Then, $\min \eqref{problem:revenue_std} = \max \eqref{problem:revenue_sdp1}$. 
\end{proposition}

\begin{proof}
    The proof is similar to Proposition \ref{prop:newsvendor_uq1}. 
\end{proof}

Similar to the Newsvendor application, higher-order moments can be incorporated, and the resulting SDP reformulation can be obtained. 

The maximum revenue is equal to $[-\max \eqref{problem:revenue_sdp1}]$. When the first and second-order moments of the supply quantity are at most $2$, i.e., $\gamma_1 = \gamma_2 = 2$, the maximum expected revenue is \$6.6495. 

We recover an optimal measure via the SDP below, \begin{align} \label{problem:revenue_recovery}
    \min_{\substack{y^k \in \R^5, z^k \in \R \\ k=1,\ldots,2n}}\ & \sum_{k=1}^{2n} z^k \\ \st\quad\ & \sum_{k=1}^{2n} y^k_1 \leq \gamma_1, \ \sum_{k=1}^{2n} y_2^k \leq \gamma_2,\ \sum_{k=1}^{2n} y_0^k = 1, \nonumber \\ & [\alpha_k b_k^2 + \beta_k b_k^4 + c_k] y_0^k  - [ 2\alpha_k b_k + 4 \beta_k b_k^3] y_1^k  +  [\alpha_k + 6 \beta_k b_k^2] y_2^k - \nonumber \\ & \quad\quad\quad 4 \beta_k b_k y_3^k + \beta_k y^k_4 \leq z^k,\ k=1,\ldots,n, \nonumber \\ & [\alpha_k b_k^2 + \beta_k b_k^4 + c_k] y_0^{n+k} - [\alpha_k b_k + \beta_k b_k^3] y_1^{n+k}  \leq z^{n+k},\ k=1,\ldots,n, \nonumber \\ & c_k y^{n+k}_0 \leq z^{n+k},\ k=1,\ldots,n, \nonumber \\ & y_0^k F_0 + y_1^k F_1 \succeq 0,\ \begin{pmatrix}
        y_0^k & y_1^k & y_2^k \\ y_1^k & y_2^k & y_3^k \\ y_2^k & y_3^k & y_4^k
    \end{pmatrix} \succeq 0,\ k=1,\ldots,2n. \nonumber
\end{align}

Solving \eqref{problem:revenue_recovery} for $n=3$ gives $\bar y^1 = (\bar y^1_0, \bar y_1^1 , \bar y_2^1 , \bar y_3^1, \bar y_4^1) = (0,0,0,0,0)$, $\bar y^2 = (1, 1.4142, 2, 2.8284, 4)$, $\bar y^3 = (0,0,0,0,0)$, $\bar y^4 = (0,0,0,0,1.5563)$, $\bar y^5 = (0,0,0,0,1.5871)$, $\bar y^6 = (0,0,0,0,1.4684)$, $\bar z^k = 0$, $k=1,3,4,5,6$, and $\bar z^2 = -6.6495$.  Clearly, $\bar y_0^k \geq 0$ and $K=\{k \in \{1,\ldots,6\} : \bar y_0^k>0\}=\{2\}$. Moreover, $\bar z^k \geq 0$ for all $k \notin K$. For $h_1 (v) = v$, it is satisfied that $\sum_{\alpha=  0}^4 \bar y_\alpha^k (h_1)_\alpha = \bar y_1^k = 0$ for all $k \notin K$, and for $h_2 (v) = v^2$, $\sum_{\alpha = 0}^4 \bar y_\alpha^k (h_2)_\alpha = \bar y_2^k = 0$ for all $k \notin K$. By Theorem \ref{thm:uq_recovery}, an optimal probability measure for problem \eqref{problem:revenue_std}, and thus for \eqref{problem:revenue_uq} is $\1_{1.4142}$, which is the discrete distribution $\mathbb{P} (\omega = 1.4142) = 1$. This result can be interpreted by the fact that the maximum expected revenue $\$6.6495$ can be achieved when the merchant supplies $1.4142$ units of goods. 

We want to remark that the assumption $\bar y_0^k > 0$, $k=1,\ldots,r$, in Corollary \ref{coro:recovery2} is sufficient but not necessary for an optimal solution to be $\sum_{k=1}^{r} \bar y_0^k \1_{\widehat u^k}$. Indeed, with $r=2n=6$, the above revenue problem offers an example that $\bar y_0^k > 0$, $k=1,\ldots,r$, is not a necessary condition. 

\begin{figure}[H]
    \centering    
    \includegraphics[width=0.4\linewidth]{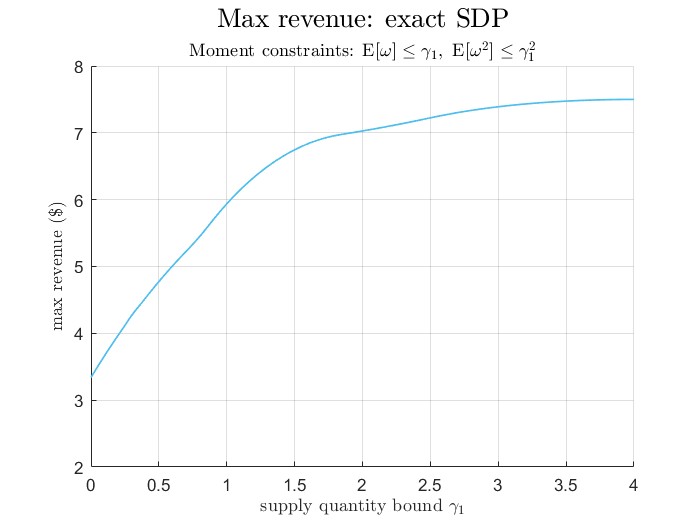}
    \includegraphics[width=0.4\linewidth]{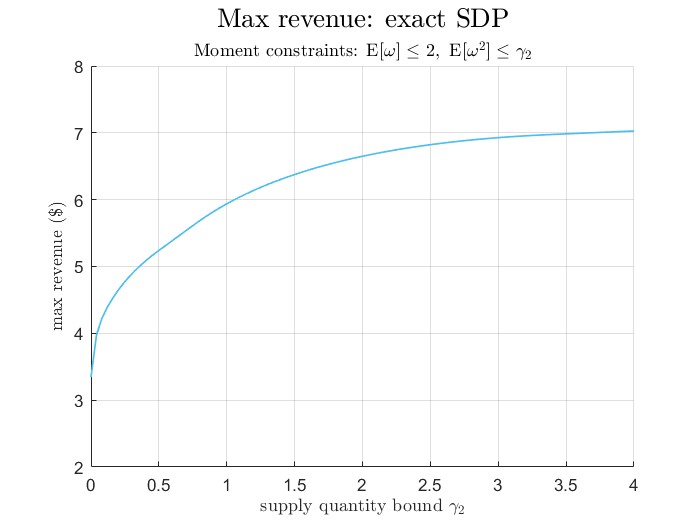}
    \caption{Maximum expected revenue by varying $\gamma_1$ (left) and $\gamma_2$ (right). }
    \label{fig:revenue-sdp}
\end{figure}

Figure \ref{fig:revenue-sdp} (left) shows how the maximum expected revenue changes with respect to the mean bound $\gamma_1$. The bound for the second-order moment is set to be $\gamma_2 = \gamma_1^2$. The curve exhibits an increasing trend. For a large mean bound $\gamma_1$, the company can target the customer offering a higher price, and thus the revenue is higher. Figure \ref{fig:revenue-sdp} (right), on the contrary, fixes $\gamma_1 = 2$ and increases $\gamma_2$. The curve exhibits an upward shape, which can be justified by the fact that the larger the second-order moment (and thus the variance), the wider the customers the company can sell products to, and thus the higher the revenue. 

\section{Conclusion and Outlook }\label{sec:conclusion}

We showed how to derive Sum-of-Squares (SOS) reformulations for important classes of infinite-dimensional moment optimization problems involving piecewise SOS-convex functions. We also showed how to recover an optimal probability measure from the associated Semi-Definite Program (SDP) reformulation. 

It is worth noting that the class of piecewise SOS-convex functions as given in Definition \ref{defn:piecewise-sos} facilitates tractable representations of the corresponding moment problem (Theorem \ref{thm:uqd_max_lin}) and it is rich enough to cover a broad class of functions encountered in optimization. Our definition \ref{defn:piecewise-sos} of piecewise functions may be extended to more general settings where one piece of SOS-convex polynomial is given at one partition of $\R^m$, in line with the piecewise linear-quadratic function studied in \cite[Definition~10.20]{rockafellar2009variational}. It would be intriguing to explore the applicability and tractability of these generalized classes of piecewise functions.

Further, our method in Theorem \ref{thm:max_sos_conv_spect} suggests that the numerically tractable SOS representation for non-negative piecewise SOS-convex functions over a projected spectrahedron can be extended to more general settings. For instance, it would be of interest to examine similar representations for piecewise SOS-convex functions over any non-convex sets whose convex hulls are semi-definite-representable, leading to SDP reformulations or relaxations for broad classes of optimization problems. 

As applications, we presented numerical results for a class of generalized Newsvendor and revenue maximization problems with higher-order moments by solving their equivalent SDP reformulations. Our approach opens new avenues for further research, such as conic program reformulations for distributionally robust optimization problems \cite{wiesemann2014distributionally,zhen2023unified} involving piecewise SOS-convex functions, and applications to practical models such as the lot-sizing and product management problems in the face of uncertain conditions \cite{ben2009robust}. These problems will be examined in our forthcoming studies. 

\section*{Acknowledgment}
The authors are grateful to the referees for their helpful comments and valuable suggestions which have contributed to the final preparation of the paper.

\section*{Declarations}

\textbf{Data availability statement}. No data was used for the research described in the article. 

\textbf{Conflict of interest statement}. The authors have no conflict of interest to declare that are relevant to the content of this article.

\begin{appendices}
    \section{Piecewise SOS-Convex Representation } \label{sec:piecewise}

In this appendix, we present explicit piecewise SOS-convex representations for the examples mentioned in the Introduction (Section \ref{sec:intro}). We note that the representation is, in general, not unique. \begin{itemize}

    \item The truncated $\ell_1$-norm, $p_\varepsilon (v) = \min \{1 , \varepsilon |v| \}$, $\varepsilon > 0$, can be expressed as $\displaystyle p_\varepsilon (v) = \min_{k=1,2} \max_{\ell=1,2} g_\ell^k (v)$ where $g_1^1 (v) = g_2^1 (v) = 1, g_1^2 (v) = \varepsilon v, g_2^2(v) = -\varepsilon v$. Note that the truncated $\ell_1$-norm is non-convex and non-smooth. We plot $\displaystyle \max_{\ell=1,2} g_\ell^1 (v)$ versus $\displaystyle \max_{\ell=1,2} g_\ell^2 (v)$ in Figure \ref{fig:lasso} (left) and the truncated $\ell_1$-norm in Figure \ref{fig:lasso} (right). 

    \begin{figure}[H]
        \centering
        \includegraphics[width=0.4\linewidth]{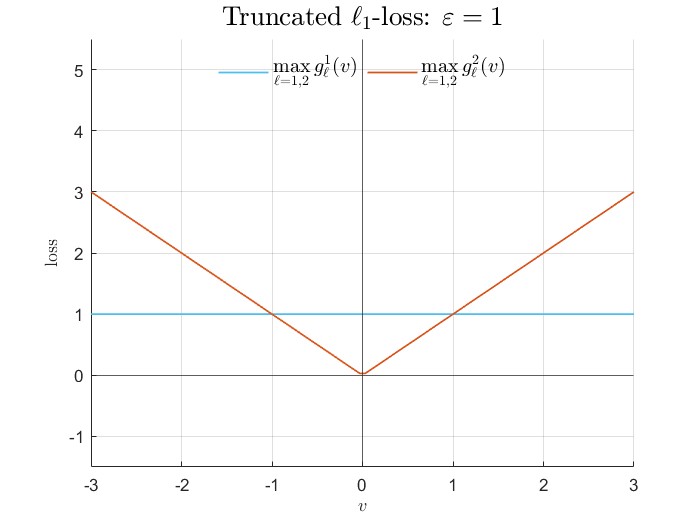}
        \includegraphics[width=0.4\linewidth]{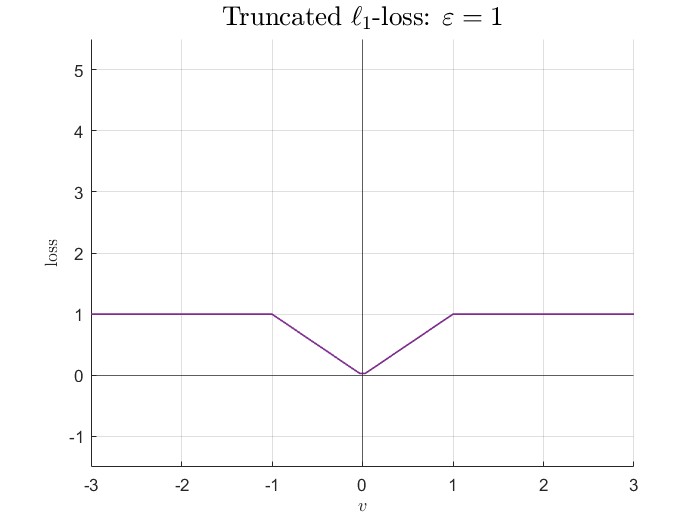}
        \caption{Piecewise representation (left) and truncated $\ell_1$-norm (right). }
        \label{fig:lasso}
    \end{figure}
    
    \item The piecewise quadratic function, $a,b \geq 0$, $c\in \R$,  \begin{align*}
        f (v) = \begin{cases}
            a (v-b)^2 + c,\ & \text{if}\ 0 \leq v \leq b, \\ c,\ & \text{if}\ v > b,
        \end{cases}
    \end{align*}can be written as $\displaystyle f(v) = \min_{k=1,2} \max_{\ell=1,2} g_\ell^k (v)$ where $g_1^1( v) = g_2^1 (v) = a (v-b)^2 + c$, $g_1^2(v) = -ab v+ab^2+c$, and $g_2^2(v) = c$. We plot $\displaystyle \max_{\ell=1,2} g_\ell^1 (v)$ versus $\displaystyle \max_{\ell=1,2} g_\ell^2 (v)$ in Figure \ref{fig:disutility} (left) and the full piecewise function in Figure \ref{fig:disutility} (right). 

    \begin{figure}[H]
        \centering
        \includegraphics[width=0.4\linewidth]{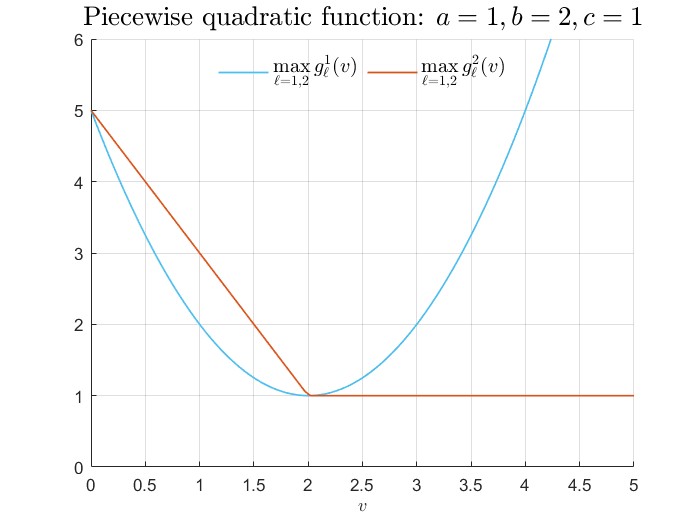}
        \includegraphics[width=0.4\linewidth]{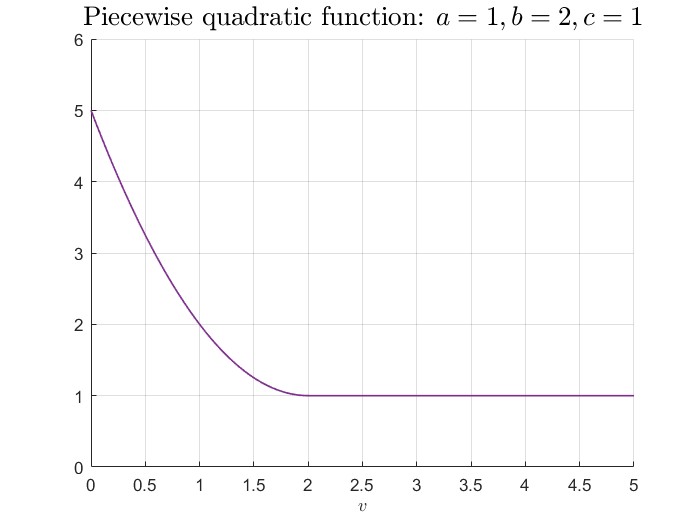}
        \caption{Piecewise representation (left) and the piecewise quadratic function (right). }
        \label{fig:disutility}
    \end{figure}

    \item The piecewise quadratic-quartic function from Section \ref{sec:applications}, \begin{align*}
        f(v) = \begin{cases}
            \alpha (v-b)^2 + \beta (v-b)^4 + c,\ & \text{if}\ 0 \leq v \leq b, \\ c , \ & \text{if}\ v > b,
        \end{cases}
    \end{align*}with $\alpha, \beta, b \geq 0$, $c\in \R$, can be expressed as $\displaystyle f(v) = \min_{k=1,2} \max_{\ell=1,2} g_\ell^k (v)$, where $g_1^1 (v) = g_2^1 (v) = \alpha (v- b)^2 + \beta (v-b)^4 + c$, $g_1^2 (v) = -(\alpha b + \beta b^3) v + (\alpha b^2 + \beta b^4 + c)$, and $g_2^2 (v) = c$. See Figure \ref{fig:disutility2} for illustrations. 
    \begin{figure}[H]
        \centering
        \includegraphics[width=0.4\linewidth]{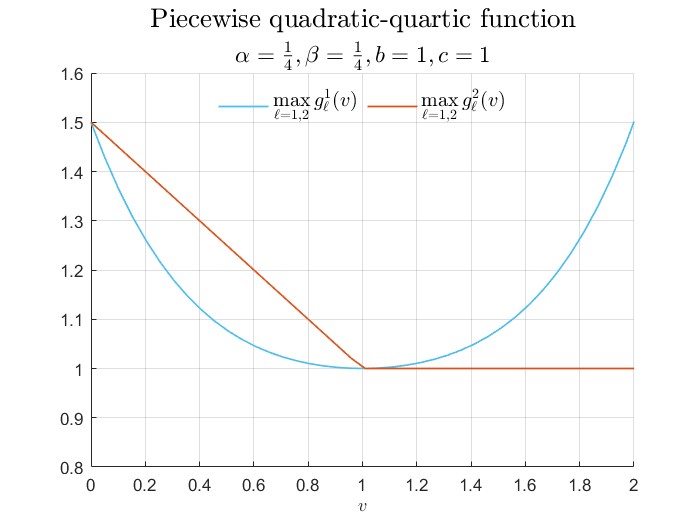}
        \includegraphics[width=0.4\linewidth]{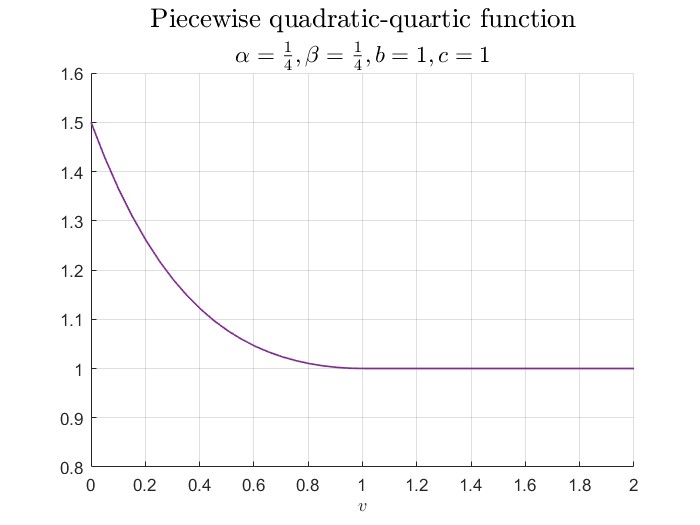}
        \caption{Piecewise representation (left) and the piecewise quadratic-quartic function (right). }
        \label{fig:disutility2}
    \end{figure}

    \item The Huber loss function \cite{wiesemann2014distributionally}  with parameter $\varepsilon > 0$, \begin{align*}
        H_\varepsilon (v) = \begin{cases}
            \frac{1}{2} v^2,\ & \text{if}\ |v| \leq \varepsilon, \\ \varepsilon |v| - \frac{1}{2} \varepsilon^2,\ & \text{otherwise}, 
        \end{cases}
    \end{align*}can be expressed equivalently as $\displaystyle H_\varepsilon(v) = \min_{k=1,2} \max_{\ell=1,2,3,4} g_\ell^k (v)$, where $g_\ell^1 (v) = \frac{1}{2} v^2$, $\ell=1,\ldots,4$, $g_1^2(v) = \varepsilon v - \frac{1}{2} \varepsilon^2$, $g_2^2(v) = -\varepsilon v - \frac{1}{2} \varepsilon^2$, $g_3^2 (v) = \frac{1}{2} \varepsilon v$, and $g_4^2 (v) = -\frac{1}{2} \varepsilon v$. See Figure \ref{fig:huber} for illustrations. 

    \begin{figure}[H]
        \centering
        \includegraphics[width=0.4\linewidth]{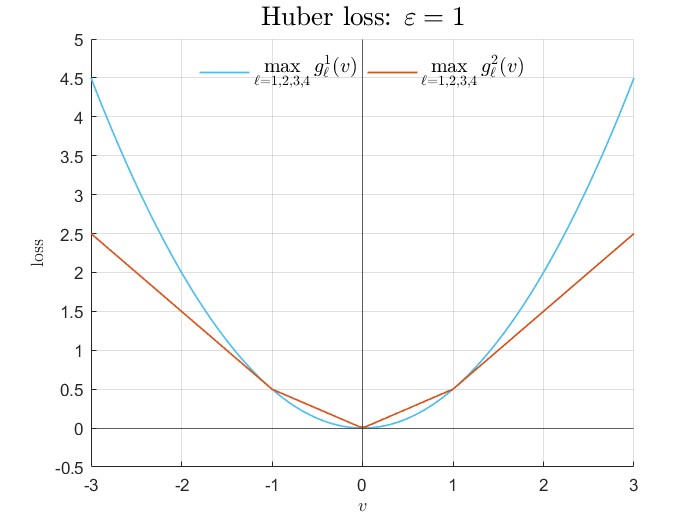}
        \includegraphics[width=0.4\linewidth]{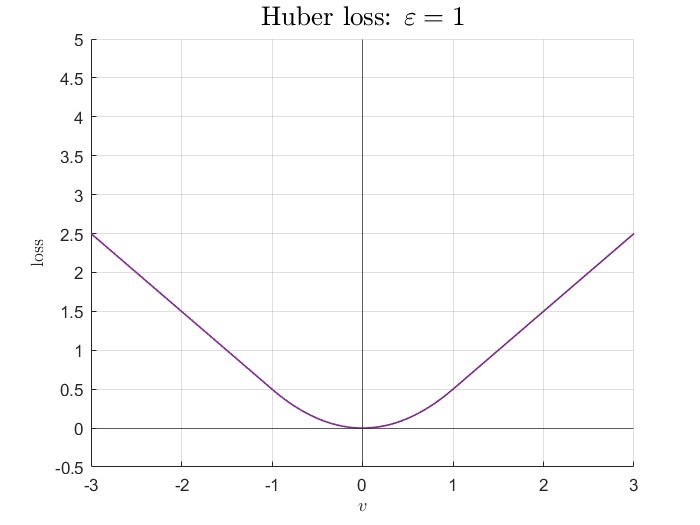}
        \caption{Piecewise representation (left) and Huber loss function (right). }
        \label{fig:huber}
    \end{figure}
\end{itemize}

\section{Conic Duality in Topological Vector Spaces} \label{sec:conic}

This section presents duality results for infinite-dimensional conic linear programs that were used in Section \ref{sec:mmt}. 

Let $X$ and $X'$ be real topological vector spaces. They are paired if a bilinear form $\langle \cdot, \cdot \rangle : X' \times X \to \R$ is defined. For the linear map $A : X \to \R^J$, assume that for any $\lambda \in \R^J$, there is a unique $x' \in X'$ satisfying $\lambda^\top A x =\langle x', x\rangle$ for all $x \in X$. See \cite[Assumption~A1]{shapiro2001duality} for a discussion of this assumption. For the continuous linear map $A : X \to \R^J$, the adjoint mapping $A^* : \R^J \to X'$ is defined as $\langle A^*\lambda, x \rangle = \lambda^\top A x$ for any $x \in X, \lambda\in \R^J$. See \cite{bonnans2013perturbation,zalinescu2002convex} for more details on the convexity of sets and functions. 

We consider the problem \begin{align} \label{problem:conic} \tag{CP}
    \min_{x \in X} \ & \langle v,x\rangle \\ \st \ & x \in C, -Ax + b \in K, \nonumber
\end{align}where $C\subset X$ and $K \subset \R^J$ are convex cones that are closed in the respective topologies, $A : X \to \R^J$ is a linear map, and $v\in X', b \in \R^J$. 

Associate \eqref{problem:conic} with the set $D = \{ x \in X : -Ax + b\in K\}$ and the dual problem \begin{align} \label{problem:conic_dual} \tag{CD} 
    \max_{\lambda\in \R^J}\ & -\lambda^\top b \\ \st\ \ & v + A^*\lambda\in C^+, \lambda\in K^+, \nonumber
\end{align}where $C^+=\{x' \in X' : \langle x', x\rangle \geq 0, \forall x\in C\}$ and $K^+=\{w \in \R^J: w ^\top x \ge 0, \forall x \in K\}$. 

We present a strong duality theorem under a more general constraint qualification known as the Generalized-Sharpened strong Conical Hull Intersection Property (G-S strong CHIP \cite{chieu2018constraint}). A pair of sets $\{C, D\}$ is said to satisfy the G-S strong CHIP at $x^*\in ( C\cap D)$ whenever 
$(C\cap D - x^*)^+ = (C-  x^*)^+ + \cup_{\lambda\in K^+} \{- A^* \lambda : \lambda^\top (-A x^* + b) = 0\}$.

\begin{theorem}[{\bf Strong duality under G-S strong CHIP}] \label{thm:conic_dual_chip}
    Suppose that $x^*$ is a minimizer for the problem \eqref{problem:conic}, the linear form $\langle u,\cdot \rangle$ is continuous at $x^*$ for each $u\in X'$, the linear map $A : X\to \R^J$ is continuous, the convex cones $C$ and $K$ are closed in the respective topologies, and that 
    the pair $\{C, D\}$ satisfies the G-S strong CHIP at $x^*$. Then, $\min \eqref{problem:conic} = \max \eqref{problem:conic_dual}$, and the maximum of $ \eqref{problem:conic_dual}$ is attained. 
\end{theorem}

\begin{proof} 
Firstly, we note that, by construction, weak duality holds, that is, $\min \eqref{problem:conic} \geq \max \eqref{problem:conic_dual}$. To see $\min \eqref{problem:conic} \leq \max \eqref{problem:conic_dual}$, let $x^*$ be a minimizer of \eqref{problem:conic}. By the necessary optimality conditions \cite{zalinescu2002convex}, $0 \in  \partial (\langle v,\cdot \rangle)(x^*) - (C\cap D - x^*)^+= \{v\} - (C\cap D - x^*)^+$. This means $v \in (C\cap D - x^*)^+$, and thus $v\in (C - x^*)^+ + \cup_{\lambda\in K^+} \{- A^* \lambda : \lambda^\top (-Ax^* + b) = 0\}$ by the G-S strong CHIP assumption. Therefore, there exist $u \in (C-x^*)^+$ and $\bar {\lambda}\in K^+$ satisfying $u = v + A^* \bar {\lambda} \in (C-x^*)^+$, $\bar {\lambda}^\top (-A x^* + b) = 0$, and  \begin{equation}\label{eq:use0}         \langle v + A^* \bar {\lambda}, x - x^*\rangle \geq 0,\ \mbox{ for all }\ x \in C.\end{equation} 

Since $C$ is a convex cone, we have $\langle v + A^* \bar \lambda, x\rangle \geq 0$ for all $x\in C$, and so $v + A^* \bar \lambda \in C^+$. In particular, $\bar {\lambda}$ is feasible for \eqref{problem:conic_dual}. Moreover, letting $x=0\in C$ in Eqn. \eqref{eq:use0} gives $\langle v + A^*\bar {\lambda}, x^* \rangle = \langle v,x^*\rangle + \bar {\lambda}^\top A x^* \leq 0$. But $\bar {\lambda}^\top Ax^* = \bar {\lambda}^\top b$, which implies $\langle v,x^*\rangle + \bar {\lambda}^\top b\leq 0$. Hence, $- \bar {\lambda}^\top b\geq \langle v,x^* \rangle$, and so,  $$\max \eqref{problem:conic_dual}=\max \{-\lambda^\top b : \lambda \in K^+, \, v + A^*\lambda \in C^+\} \ge -\bar {\lambda}^\top b \geq \langle v,x^* \rangle = \min \eqref{problem:conic}.$$ 
Therefore, we see that $\min\eqref{problem:conic}=\max \eqref{problem:conic_dual}$ and the maximum of $\eqref{problem:conic_dual}$ is attained at $\bar {\lambda}$.
\end{proof}

The G-S strong CHIP is the weakest constraint qualification guaranteeing strong (Lagrangian) duality for convex optimization. For details, see \cite{chieu2018constraint} and other references therein. As we see below, 
strong duality under the interior point condition \cite{shapiro2001duality} is a consequence of Theorem \ref{thm:conic_dual_chip}.

\begin{corollary}[{\bf Strong duality under interior point condition}] \label{coro:int_pt} 
    Suppose that $x^*$ is a minimizer for the problem \eqref{problem:conic}, the linear form $\langle u,\cdot \rangle$ is continuous at $x^*$ for each $u\in X'$, the linear map $A : X \to \R^J$ is continuous, the convex cones $C$ and $K$ are closed in the respective topologies. If $b \in \interior (A(C) + K)$, then, $\min\eqref{problem:conic} = \max\eqref{problem:conic_dual}$, and the maximum of $\eqref{problem:conic_dual}$ is attained.    
\end{corollary}

\begin{proof} 
    The conclusion will follow from Theorem \ref{thm:conic_dual_chip} if we show that $b\in \interior (A(C) + K)$ implies \begin{align*}
        (C \cap D - x^*)^+ = (C - x^*)^+ + \cup_{\lambda\in K^+} \{-A^* \lambda : \lambda^\top (-A x^* + b) = 0\}.
    \end{align*}
    
    Note that the inclusion $(C - x^*)^+ + \cup_{\lambda \in K^+} \{-A^* \lambda: \lambda^\top (-A x^* + b) = 0\} \subseteq (C\cap D - x^*)^+$ holds by construction. Conversely, let $u \in (C\cap D - x^*)^+$. By definition, this means $\langle u,x-x^* \rangle \geq 0$ for all $x\in (C\cap D)$. Hence, $\langle u,x\rangle \geq \langle u,x^* \rangle$ for all $x\in (C\cap D)$, and $x^*$ is a minimizer of the convex optimization problem $\min \{\langle u,x\rangle : x\in C, -A x + b\in K\}$. By \cite{shapiro2001duality}, there exists $\lambda \in K^+$ with $\langle u,x^* \rangle = -\lambda^\top b$ and $A^* \lambda + u \in C^+$. Since $- Ax^* + b\in K$ and $\lambda \in K^+$, we have $\lambda^\top (-A x^* + b) \geq 0$. 
    
    On the other hand, $\lambda^\top (-A x^* + b) = \langle -A^* \lambda - u, x^* \rangle \leq 0$ as $x^* \in C$ and $A^* \lambda + u\in C^+$. These together force $\lambda^\top (-Ax^* + b) = 0$. In addition, $\langle A^*\lambda+ u, x- x^*\rangle = \langle A^* \lambda + u,x\rangle \geq 0$ for all $x\in C$. This is equivalent to $\langle A^* \lambda + u, x\rangle \geq 0$ for all $x\in (C - x^*)$, which further implies $A^* \lambda + u \in (C-x^*)^+$. Hence, $u \in (C-x^*)^+ + \cup_{\lambda\in K^+} \{ -A^* \lambda : \lambda^\top (-Ax^* + b) = 0\}$, and the proof is complete.     
\end{proof}

\section{SDP Representation of SOS Problems } \label{sec:sos-sdp}

It is known that the SOS optimization problem \eqref{problem:uq_dual} can be expressed equivalently as an SDP. This appendix provides technical results for Section \ref{sec:mmt}. 

A monomial over $v \in \R^m$ of degree $\bar d$ is $v^\alpha = v_1^{\alpha_1} v_2 ^{\alpha_2} \ldots v_m^{\alpha_m}$ with $\bar d = \sum_{i=1}^{m} \alpha_i$, and $\alpha = (\alpha_1, \ldots, \alpha_m) \in (\{0\} \cup \mathbb N)^{m}$ is a multi-index. The canonical basis is \begin{align*}
    y(v) := (1, v_1,\ldots,v_m, v_1^2, v_1 v_2, \ldots, v_m^2, \ldots, v_1^{\bar d}, \ldots, v_m^{\bar d})^\top,
\end{align*}which is of dimension $s(m, \bar d) := \binom{m+\bar d}{\bar d}$. 

Let $f$ be a real polynomial with an even degree $d = 2\bar d$ written as $f (v) = \sum_{\alpha \in \mathcal{N}} (f)_\alpha v^\alpha$, where $(f)_\alpha$ is the $\alpha$-th coefficient of $f$ and  $\mathcal{N} = \{\alpha \in (\{0 \} \cup \mathbb N)^{m} : \sum_{i=1}^m \alpha_i \leq d\}$. By \cite[Proposition~2.1]{lasserre2009moments}, $f$ is SOS if and only if there exists $Q \in \mathbb S_+^{s(m,d/2)}$ such that $f (v) = y(v)^\top Q y(v)$. This expresses the coefficients of $f (v)$ as linear equations of the entries in $Q$. If we write $y(v) y(v)^\top = \sum_{\alpha \in \mathcal{N}} B_\alpha v^\alpha$ for appropriate matrices $B_\alpha \in \mathbb S^{s(m,d/2)}, \alpha \in \mathcal{N}$, checking whether $f$ is SOS amounts to finding $Q \in \mathbb S_+^{s(m,d/2)}$ such that $\tr(QB_\alpha) = (f)_\alpha$ for all $\alpha \in \mathcal{N}$. 

The following proposition shows that the SOS problem \eqref{problem:uq_dual} and the SDP program \eqref{problem:uq_sdp} share the same optimal values.

\begin{proposition} \label{prop:sos-sdp} 
    Let $g_\ell^k$, $\ell=1,\ldots,L$, $k=1,\ldots,r$, $h_j$, $j=1,\ldots,J$, be defined as in Theorem \ref{thm:uqd_max_lin}. Then, $\max \eqref{problem:uq_dual} = \max \eqref{problem:uq_sdp}$. 
\end{proposition}

\begin{proof}
The SOS constraints \begin{align*}
         \sum_{\ell=1}^{L} \delta_\ell^k  g_\ell^k (v) + \sum_{j=1}^{J} \lambda_{j} h_j (v) + \lambda_{J+1}  - \tr(Z_kF_0) -  \sum_{i=1}^m v_i \tr(Z_kF_i) \in \Sigma_d^2 (v),\ k=1,\ldots,r,
     \end{align*}of \eqref{problem:uq_dual} are equivalent to the existence of $Q_k \in \mathbb S_+^{s(m,d/2)}$, $k=1,\ldots,r$, such that  \begin{align*}
         \sum_{\ell=1}^L \delta_\ell^k (g_\ell^k)_\alpha + \sum_{j=1}^J \lambda_j (h_j)_\alpha + \lambda_{J+1} (1)_\alpha - \tr(Z_k F_\alpha) = \tr(Q_k B_\alpha),\  \text{for all}\ \alpha \in \mathcal{N},\ k=1,\ldots,r,
     \end{align*}where $F_0 := F_0 \in \mathbb{S}^{\nu}$, $F_{e_i} := F_i \in \mathbb{S}^{\nu}$, $i=1,\ldots,m$, and $F_\alpha$, $\alpha \in \mathcal{N} \setminus \{0, e_1,\ldots,e_m\}$, is the zero matrix. Thus, the conclusion follows. 
\end{proof}

\end{appendices}
\footnotesize{
\bibliography{mybib}
\bibliographystyle{abbrv}}

\end{document}